\documentclass[11pt]{article}
\usepackage[a4paper]{geometry}
\usepackage{amssymb,latexsym,amsmath,amsfonts,amsthm}
\usepackage{graphicx,epstopdf}
\usepackage{dsfont}
\usepackage{mathrsfs}
\usepackage{enumerate}
\usepackage{epsfig}
\usepackage{booktabs}
\usepackage{mathtools}
\usepackage{multirow}
\usepackage{pifont}
\usepackage{comment,verbatim}
\usepackage{overpic}
\usepackage{hyperref}
\usepackage{bbm}
\usepackage[toc,page]{appendix}

\newtheorem{theorem}{Theorem}[section]
\newtheorem{lemma}[theorem]{Lemma}

\newtheorem{corollary}[theorem]{Corollary}

\theoremstyle{definition}

\theoremstyle{remark}

\newtheorem{remark}[theorem]{Remark}

\numberwithin{equation}{section}

\usepackage{color}

\usepackage[normalem]{ulem}

\begin{document}
\title{Analysis of error localization of Chebyshev spectral approximations}
\author{Haiyong Wang\footnotemark[1]~\footnotemark[2]}
\date{}
\maketitle
\renewcommand{\thefootnote}{\fnsymbol{footnote}}
\footnotetext[1]{School of Mathematics and Statistics, Huazhong
University of Science and Technology, Wuhan 430074, P. R. China.
E-mail: \texttt{haiyongwang@hust.edu.cn}}

\footnotetext[2]{Hubei Key Laboratory of Engineering Modeling and
Scientific Computing, Huazhong University of Science and Technology,
Wuhan 430074, P. R. China.}

\begin{abstract}
Chebyshev spectral methods are widely used in numerical
computations. When the underlying function has a singularity, it has
been observed by L. N. Trefethen in 2011 that its Chebyshev
interpolants exhibit an error localization property, that is, their
errors in a neighborhood of the singularity are obviously larger
than elsewhere. In this paper, we first present a pointwise error
analysis for Chebyshev projections of functions with a singularity
and prove that the rate of convergence of Chebyshev projections of
degree $n$ at each point away from the singularity is one power of
$n$ faster than that of at the singularity. This gives a rigorous
justification for the error localization of Chebyshev projections.
We then extend the framework of our analysis to Chebyshev
interpolants, Chebyshev spectral differentiations and Legendre
projections and justify their error localization using similar
arguments. As a result, we find that Chebyshev spectral
differentiations converge faster than their best counterparts except
in a neighborhood of the singularity and, in the particular case
where the singularity is located in the interior of interval, they
converge even faster than their best counterparts in the maximum
norm.
\end{abstract}


\vspace{0.05in}


\section{Introduction}\label{sec:introduction}
Spectral approximations, such as Chebyshev and Legendre projections
and interpolants, are invaluable and powerful methods and they play
an important role in numerous practical applications, including
Gauss and Clenshaw-Curtis quadrature, rootfinding and spectral
methods for differential and integral equations (see, e.g.,
\cite{boyd2000cheb,boyd2013zero,canuto2006spect,mason2003cheb,rivlin2020cheb,shen2011spect,trefethen2000spect,trefethen2008clen,trefethen2020atap}).
Of particular importance are Chebyshev spectral approximations,
which have several remarkable advantages including: (i) they are
near-best approximations, that is, their approximation quality in
the maximum norm is close to that of best approximations
\cite{li2004optimal,mason2003cheb,rivlin2020cheb,trefethen2020atap,wang2020legendre};
(ii) the discrete Chebyshev transforms, i.e., the transforms between
the values at the Chebyshev points and the Chebyshev expansion
coefficients, can be achieved rapidly by means of the fast Fourier
transform (FFT) \cite{mason2003cheb,shen2011spect}; (iii) the
evaluation of Chebyshev interpolants can be achieved stably and
rapidly using the barycentric formula
\cite{berrut2004bary,higham2004bary}. Due to these attractive
advantages, Chebyshev spectral approximations are widely used in
many branches of numerical analysis. We refer the interested reader
to
\cite{driscoll2014cheb,mason2003cheb,rivlin2020cheb,trefethen2020atap}
for more extensive overview.

Let $\mathrm{\Omega}:=[-1,1]$ and let $T_k(x)$ be the Chebyshev
polynomial of the first kind of degree $k$, i.e.,
$T_k(x)=\cos(k\arccos(x))$. It is well known that the Chebyshev
polynomials $\{T_k(x)\}$ are orthogonal with respect to the weight
function $(1-x^2)^{-1/2}$ on $\mathrm{\Omega}$. If $f$ satisfies the
Dini-Lipschitz continuous on $\mathrm{\Omega}$, then it has the
following uniformly convergent Chebyshev series
\cite[Theorem~5.7]{mason2003cheb}
\begin{align}\label{def:ChebExp}
f(x) = \sum_{k=0}^{\infty}{'} a_k T_k(x), \quad a_k = \frac{2}{\pi}
\int_{\Omega} \frac{f(x) T_k(x)}{\sqrt{1-x^2}} \mathrm{d}x,
\end{align}
where the prime indicates that the first term of the sum should be
halved. Truncating the above infinite series after the first $n+1$
terms, we obtain the Chebyshev projection of degree $n$ of $f$,
i.e.,
\begin{align}\label{eq:ChebProj}
f_n(x) = \sum_{k=0}^{n}{'} a_k T_k(x).
\end{align}
Let $\Pi_n$ denote the space of polynomials of degree at most $n$
and let $p_n^{*}$ denote the best approximation in $\Pi_n$ to $f$ in
the maximum norm, i.e., $\|f-p_n^{*}\|_{L^{\infty}(\Omega)} =
\min_{h\in\Pi_n} \|f - h\|_{L^{\infty}(\Omega)}$. It is well known
that $p_n^{*}$ exists and is unique whenever $f\in
\mathrm{C}(\mathrm{\Omega})$. We are interested in the comparison of
$f_n$ and $p_n^{*}$. From the viewpoint of minimizing the maximum
error of approximants, it is evident that $p_n^{*}$ is better than
$f_n$. From the viewpoint of practical convenience, however, it is
evident that $f_n$ is preferable to $p_n^{*}$ since $f_n$ depends
linearly on $f$ and its implementations can be achieved efficiently
by means of the FFT. In contrast, $p_n^{*}$ depends nonlinearly on
$f$ and its implementation must resort to iterative methods, which
results in an expensive computational cost, especially when $n$ is
large. We further consider the difference in approximation quality
between $p_n^{*}$ and $f_n$. A classical result associated with this
issue are the following inequalities
\cite[Theorem~3.3]{rivlin2020cheb}
\begin{align}\label{eq:Lebesgue}
\|f - p_n^{*}\|_{L^{\infty}(\Omega)} \leq \|f -
f_n\|_{L^{\infty}(\Omega)} < \left(4 + \frac{4}{\pi^2} \log n
\right) \|f - p_n^{*}\|_{L^{\infty}(\Omega)},
\end{align}
from which we see that the maximum error of $f_n$ is worse than that
of $p_n^{*}$ by at most a logarithmic factor. Direct calculation
shows that $\|f - f_n\|_{L^{\infty}(\Omega)} < 9.6 \|f -
p_n^{*}\|_{L^{\infty}(\Omega)}$ for $n\leq10^6$, which means that
the maximum error of $p_n^{*}$ will never be better than that of
$f_n$ by one digit even if the degree of both methods is one
million. Regarding \eqref{eq:Lebesgue}, we further ask the following
question: Is $p_n^{*}$ really better than $f_n$ by a logarithmic
factor? At first glance, this question has already been answered by
the inequalities above. However, they only answer this question from
the viewpoint of measuring the maximum error of both approximants.
If we compare both approximants from the viewpoint of the rate of
pointwise convergence, the inequalities above might be one-sided or
even misleading. To illustrate this, we plot the error curves of
$f_n$ and $p_n^{*}$ for the test function $f(x)=|x-1/2|^{\alpha}$
and we test two different values of $\alpha$ and $n$, respectively.
Throughout the paper, $p_n^{*}$ is calculated by the
barycentric-Remez algorithm in Chebfun (see
\cite{driscoll2014cheb}). The pointwise error curves are depicted in
Figure \ref{fig:Exam}, which suggest that:
\begin{figure}[ht]
\centering
\includegraphics[width=0.5\textwidth]{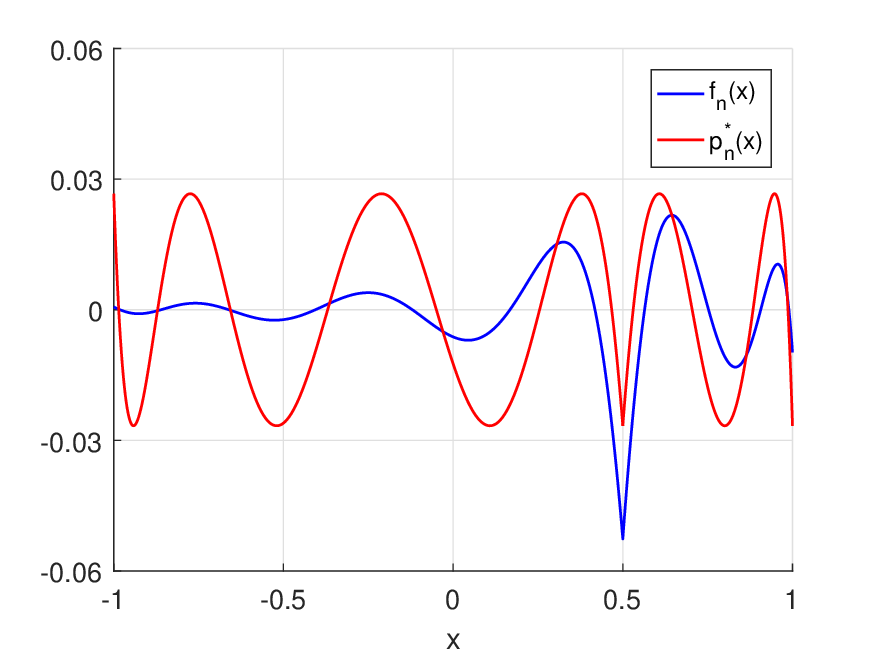}~
\includegraphics[width=0.5\textwidth]{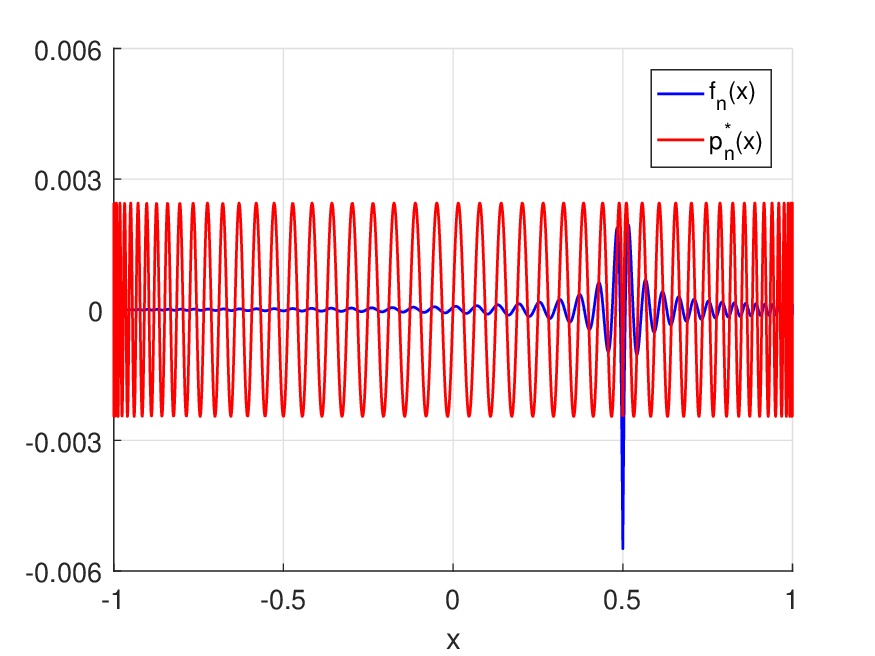}\\
\includegraphics[width=0.5\textwidth]{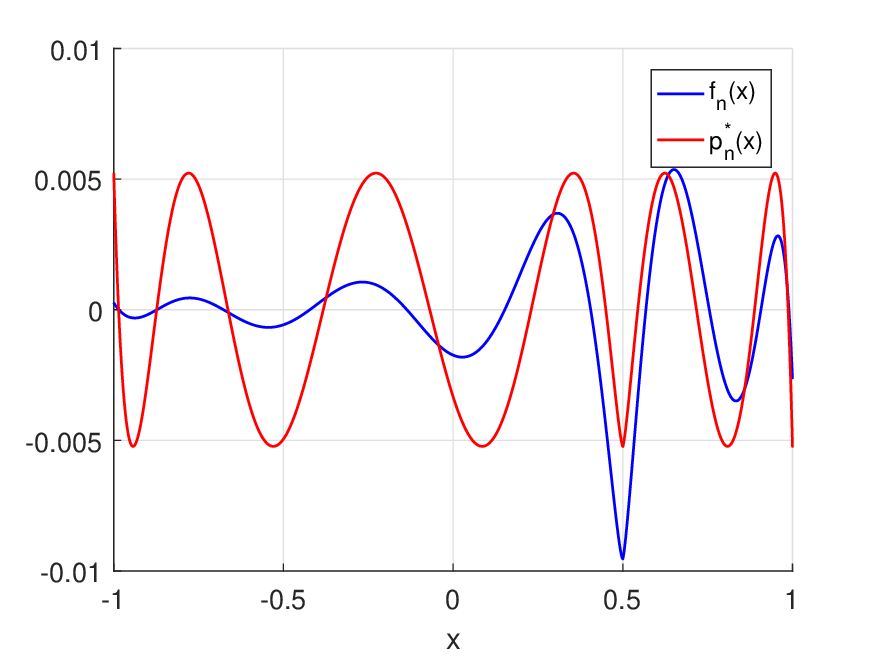}~
\includegraphics[width=0.5\textwidth]{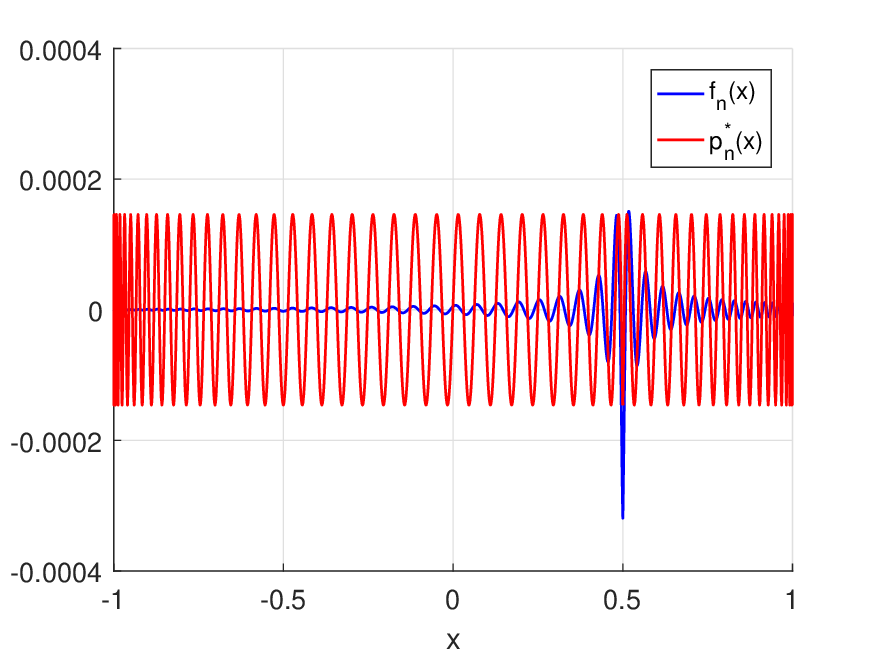}
\caption{Pointwise error curves of $f_n$ and $p_n^{*}$ for $n=10$
(left) and $n=100$ (right). The top row corresponds to $\alpha=1$
and the bottom row corresponds to $\alpha=3/2$.} \label{fig:Exam}
\end{figure}
\begin{itemize}
\item[\rm (i)] The error of $p_n^{*}$ equioscillates over the whole
interval. The error of $f_n$, however, is highly localized for large
$n$, that is, the maximum errors of $f_n$ are always attained at a
small neighborhood of the singularity and the larger $n$, the
narrower the neighborhood that the pointwise error curves attain
their maximum;

\item[\rm (ii)] Comparing the accuracy of $f_n$ and $p_n^{*}$,
we observe that the accuracy of $p_n^{*}$ is better than that of
$f_n$ only in a small neighborhood of the singularity. Otherwise,
the accuracy of $p_n^{*}$ is worse than that of $f_n$.
\end{itemize}
We are not the first to notice these observations. In fact,
Trefethen in \cite{trefethen2011six} observed a similar phenomenon
from the pointwise error curve of Chebyshev interpolant. He further
made the following comments:
\begin{quote}
``Which approximation would be more useful in an application? I
think the only reasonable answer is, it depends. Sometimes one
really does need a guarantee about worst-case behavior. In other
situations, it would be wasteful to sacrifice so much accuracy over
95\% of the range just to gain one bit of accuracy in a small
subinterval.''
\end{quote}
However, no theoretical justification was given to explain these
observations.

In this work, we aim to justify the error localization property of
Chebyshev spectral approximations and present a detailed analysis of
their pointwise rate of convergence for functions with
singularities. For ease of presentation, we restrict our analysis to
the following model function
\begin{align}\label{def:Model}
f(x)=|x-\xi|^{\alpha} g(x),
\end{align}
where $\xi\in\Omega$ and $\alpha>0$ is not an even integer whenever
$\xi\in(-1,1)$ and is not an integer whenever $\xi=\pm1$, and
$g\in{C}^{\infty}(\Omega)$ and $g(\xi)\neq0$. In order to compare
the pointwise rates of convergence of $f_n$ and $p_n^{*}$, we define
the set
\begin{align}\label{def:V}
V_{0}(x) = \left\{x\in\Omega: ~ \lim_{n\rightarrow\infty}
\frac{|f(x)-f_n(x)|}{\|f - p_n^{*}\|_{L^{\infty}(\Omega)}} = 0
\right\}.
\end{align}
Note that we have used $\|f - p_n^{*}\|_{L^{\infty}(\Omega)}$
instead of $|f(x) - p_n^{*}(x)|$ in the denominator due to the
equioscillation property of the error of $p_n^{*}$. It is easily
seen that the set $V_{0}(x)$ contains all the points on $\Omega$ at
which $f_n$ converges faster than $p_n^{*}$. Our main result states
that the rate of convergence of $f_n$ is one power of $n$ faster
than that of $p_n^{*}$ whenever $x$ is away from $\xi$ and both
$f_n$ and $p_n^{*}$ converge at the same rate whenever $x=\xi$.
Therefore, we can deduce that $V_0(x)=\Omega\backslash\{\xi\}$ and
this gives a rigorous justification of those observations displayed
in Figure \ref{fig:Exam}. We show that the key ingredient for
explaining the error localization of $f_n$ is to understand the
asymptotic behavior of the following two functions
\begin{align}\label{def:PsiFun}
\Psi_{\nu}^{\mathrm{C}}(\theta,n) = \sum_{k=n+1}^{\infty}
\frac{\cos(k\theta)}{k^{\nu+1}}, \qquad
\Psi_{\nu}^{\mathrm{S}}(\theta,n) = \sum_{k=n+1}^{\infty}
\frac{\sin(k\theta)}{k^{\nu+1}},
\end{align}
where $\theta\in\mathbb{R}$, $n\in\mathbb{N}_0:=\mathbb{N}\cup\{0\}$
and $\nu>-1$ whenever $\theta(\mathrm{mod}~ 2\pi)\neq0$ and $\nu>0$
whenever $\theta(\mathrm{mod}~ 2\pi)=0$. As will be shown later,
these two functions are intimately related to the Lerch's
transcendent function and its special case Hurwitz zeta function
(see, e.g., \cite[Chapter~25]{olver2010nist}) and their asymptotic
behavior has a crucial distinction between $\theta(\mathrm{mod}~
2\pi)\neq0$ and $\theta(\mathrm{mod}~ 2\pi)=0$. With this key
finding in hand, we are able to clearly explain not only the error
localization of Chebyshev projections, but also the same property of
Chebyshev interpolants, Chebyshev spectral differentiations and
Legendre projections. As a result, we show that Chebyshev spectral
differentiations converge faster than their best counterparts except
in a neighborhood of $\xi$ and, in the particular case where
$\xi\in(-1,1)$, they converge even faster than their best
counterparts in the maximum norm. Our results provide further
justification for the use of Chebyshev spectral methods in numerical
computations.

The paper is organized as follows. In section \ref{sec:Pointwise},
we present a thorough analysis of the pointwise error estimates of
Chebyshev projections. In section \ref{sec:Extension}, we extend our
analysis to a more general setting, including superconvergent points
of Chebyshev projections, Chebyshev interpolants, Chebyshev spectral
differentiation and Legendre projections. In section
\ref{sec:conclusion}, we finish this paper with some concluding
remarks.

\section{Pointwise error estimates of Chebyshev projections}\label{sec:Pointwise}
In this section, we present pointwise error estimates of Chebyshev
projections for the model function defined in \eqref{def:Model}. Our
main result will provide a thorough understanding of the pointwise
error behavior of Chebyshev projections for functions with a
singularity. For notational simplicity, we introduce
\begin{align}\label{def:varphi}
\varphi_{\xi}^{\pm}(x) = \frac{\arccos(x)\pm\arccos(\xi)}{2},
\end{align}
and we will drop the argument $x$ in $\varphi_{\xi}^{\pm}(x)$
whenever there is no ambiguity.

Before stating our main result, we provide two useful lemmas that
will be used in the proof of our main result.
\begin{lemma}\label{lem:AsymChebCoeff}
Let $f$ be the function defined in \eqref{def:Model} and let $a_k$
be its $k$th Chebyshev coefficient.
\begin{itemize}
\item[\rm (i)] If $\xi\in(-1,1)$, then for $k\gg1$, we have
\begin{align}
a_k = \mathcal{I}_1(\alpha,\xi) \frac{T_k(\xi)}{k^{\alpha+1}} +
\mathcal{I}_2(\alpha,\xi) \frac{U_{k-1}(\xi)}{k^{\alpha+2}} +
O(k^{-\alpha-3}),
\end{align}
where $U_k(x)$ is the Chebyshev polynomial of the second kind of
degree $k$ and
\begin{align}
\mathcal{I}_1(\alpha,\xi) &=
-\frac{4(1-\xi^2)^{\alpha/2}\Gamma(\alpha+1)}{\pi}
\sin\left(\frac{\alpha\pi}{2}\right) g(\xi),  \nonumber \\
\mathcal{I}_2(\alpha,\xi) &= -\frac{4(1-\xi^2)^{\alpha/2}
\Gamma(\alpha+2)}{\pi} \sin\left(\frac{\alpha\pi}{2}\right) \left[
(1-\xi^2) g{'}(\xi) - \frac{\alpha\xi}{2} g(\xi) \right]. \nonumber
\end{align}

\item[\rm (ii)] If $\xi=\pm1$, then for $k\gg1$, we have
\begin{align}
a_k = \mathcal{B}(\alpha,\xi) \frac{T_k(\xi)}{k^{2\alpha+1}} +
O(k^{-2\alpha-3}), \quad \mathcal{B}(\alpha,\xi) = -
\frac{\sin(\alpha\pi)\Gamma(2\alpha+1)}{2^{\alpha-1}\pi} g(\xi).
\end{align}
\end{itemize}
\end{lemma}
\begin{proof}
As for (i), it follows from Theorem 2.1 in \cite{kzaz2000asymp}. As
for (ii), it follows from Theorem 2.2 in \cite{wang2018clenshaw}.
\end{proof}

\begin{lemma}\label{lem:DecaySeries}
Let $\Psi_{\nu}^{\mathrm{C}}(\theta,n)$ and
$\Psi_{\nu}^{\mathrm{S}}(\theta,n)$ be the functions defined in
\eqref{def:PsiFun}. For $n\gg1$, the following two statements hold:
\begin{itemize}
\item[\rm (i)] If $\theta(\mathrm{mod}~ 2\pi)\neq0$ and $\nu>-1$, we have
\begin{align}\label{eq:PsiCAsyI}
\Psi_{\nu}^{\mathrm{C}}(\theta,n) &=
-\frac{\sin((2n+1)\theta/2)}{2\sin(\theta/2)} n^{-\nu-1} +
O(n^{-\nu-2}),
\end{align}
and
\begin{align}\label{eq:PsiSAsy}
\Psi_{\nu}^{\mathrm{S}}(\theta,n) &=
\frac{\cos((2n+1)\theta/2)}{2\sin(\theta/2)} n^{-\nu-1} +
O(n^{-\nu-2}).
\end{align}

\item[\rm (ii)] If $\theta(\mathrm{mod}~2\pi)=0$ and $\nu>0$, we have
\begin{align}\label{eq:PsiCAsyII}
\Psi_{\nu}^{\mathrm{C}}(\theta,n) &= \frac{1}{\nu} n^{-\nu} -
\frac{1}{2} n^{-\nu-1} + O(n^{-\nu-2}).
\end{align}
\end{itemize}
\end{lemma}
\begin{proof}
We notice that
\begin{align}
\Psi_{\nu}^{\mathrm{C}}(\theta,n) + i
\Psi_{\nu}^{\mathrm{S}}(\theta,n) &= e^{i(n+1)\theta}
\Phi(e^{i\theta},\nu+1,n+1), \nonumber
\end{align}
where $i$ is the imaginary unit and $\Phi(z,s,a)$ is the Lerch's
transcendent function defined by (see, e.g.,
\cite[\S25.14(i)]{olver2010nist} and \cite{ferreira2004asymp})
\begin{align}
\Phi(z,s,a) = \sum_{k=0}^{\infty} \frac{z^k}{(k+a)^{s}},
\end{align}
and where $a\neq0,-1,\ldots$, $s\in\mathbb{C}$ whenever $|z|<1$ and
$\Re(s)>1$ whenever $|z|=1$. For other values of $z$, $\Phi(z,s,a)$
is defined by analytic continuation. Moreover, the Lerch's
transcendent function reduces to the Hurwitz zeta functions whenever
$z=1$, i.e.,
\begin{align}
\Phi(1,s,a) = \sum_{k=0}^{\infty} \frac{1}{(k+a)^{s}} := \zeta(s,a),
\end{align}
where $\zeta(s,a)$ is the Hurwitz zeta functions
\cite[\S25.11(i)]{olver2010nist}. As for (i), we recall from
\cite[Theorem 1]{ferreira2004asymp} that $\Phi(z,s,a) = (1-z)^{-1}
a^{-s} + O(a^{-s-1})$ for $z\in\mathbb{C}\setminus[1,\infty)$,
$s\in\mathbb{C}$ and $|a|\rightarrow\infty$, and thus
\begin{align}
\Psi_{\nu}^{\mathrm{C}}(\theta,n) +
i\Psi_{\nu}^{\mathrm{S}}(\theta,n) =
\frac{e^{i(n+1)\theta}}{1-e^{i\theta}} n^{-\nu-1} +O(n^{-\nu-2}).
\nonumber
\end{align}
The results \eqref{eq:PsiCAsyI} and \eqref{eq:PsiSAsy} follow from
the above equation by taking real and imaginary parts, respectively.
As for (ii), from \cite[Equation (25.11.3)]{olver2010nist} we know
that
$\Psi_{\nu}^{\mathrm{C}}(\theta,n)=\zeta(\nu+1,n+1)=\zeta(\nu+1,n)-n^{-\nu-1}$.
The result \eqref{eq:PsiCAsyII} then follows from the asymptotic
expansion of $\zeta(\nu,n)$ in
\cite[Equation~(25.11.43)]{olver2010nist}. This completes the proof.
\end{proof}

Several remarks are in order.
\begin{remark}
%The magnitude of the Lerch's transcendent function $\Phi(z,s,a)$ has
%a peak near $z=1$; see Figure \ref{fig:Lerch} for an illustration.
When $|z|=1$ and $a\rightarrow\infty$, we notice that $\Phi(z,s,a)$
behaves like $O(a^{-s})$ whenever $z\neq1$ and $\Re(s)\in\mathbb{C}$
and behaves like $O(a^{1-s})$ whenever $z=1$ and $\Re(s)>1$. This
explains why the asymptotic behavior of
$\Psi_{\nu}^{\mathrm{C}}(\theta,n)$ has an obvious distinction
between $\theta(\mathrm{mod}~2\pi)\neq0$ and
$\theta(\mathrm{mod}~2\pi)=0$. As will be shown below, this
distinction turns out to be the key to explain the error
localization of Chebyshev spectral approximations. We mention that
this distinction also plays a key role in analyzing the rate of
pointwise convergence of the modified Fourier expansions for
nonperiodic
functions in one or more dimensions \cite{adcock2010a,olver2012}.  % in the univariate and multivariate cases
\end{remark}

\begin{remark}\label{RK:Psi}
For any $k\in\mathbb{Z}$ and $n\gg1$, we have
\begin{align}
\Psi_{\nu}^{\mathrm{C}}(k\pi,n) = \left\{
            \begin{array}{ll}
{\displaystyle O(n^{-\nu})},    & \hbox{$k$ is even,}   \\[12pt]
{\displaystyle O(n^{-\nu-1})},  & \hbox{$k$ is odd, }
            \end{array}
            \right. \nonumber
\quad \Psi_{\nu}^{\mathrm{S}}(k\pi,n) = 0.
\end{align}
\end{remark}

Our main result is stated in the following theorem.
\begin{theorem}\label{thm:PointErrorProj}
Let $f$ be the function defined in \eqref{def:Model} and let $f_n$
be its Chebyshev projection of degree $n$. As $n\rightarrow\infty$,
the following statements are true.
\begin{itemize}
\item[\rm (i)] If $\xi\in(-1,1)$, we have for $x\neq\xi$ that
\begin{align}\label{eq:AsyCaseI}
f(x) - f_n(x) &= -\frac{\mathcal{I}_1(\alpha,\xi)}{4} \left[
U_{2n}(\cos(\varphi_{\xi}^{+})) + U_{2n}(\cos(\varphi_{\xi}^{-}))
\right] n^{-\alpha-1} + O(n^{-\alpha-2}),
\end{align}
and for $x=\xi$
\begin{align}\label{eq:AsyCaseII}
f(x) - f_n(x) = \frac{\mathcal{I}_1(\alpha,\xi)}{2\alpha}
n^{-\alpha} - \mathcal{I}_1(\alpha,\xi)\frac{U_{2n}(\xi)+1}{4}
n^{-\alpha-1} + O(n^{-\alpha-2}).
\end{align}

\item[\rm (ii)] If $\xi=1$, we have
\begin{align}\label{eq:ErrAsyBI}
f(x) - f_n(x) &= \mathcal{B}(\alpha,\xi) \left\{
            \begin{array}{ll}
{\displaystyle \frac{1}{2\alpha}n^{-2\alpha} - \frac{1}{2}n^{-2\alpha-1} + O(n^{-2\alpha-2}) },  & \hbox{$x=\xi$,}   \\[12pt]
{\displaystyle -\frac{U_{2n}(\cos(\varphi_{\xi}^{+}))}{2}
n^{-2\alpha-1} + O(n^{-2\alpha-2}) }, & \hbox{$x\neq\xi$.}
            \end{array}
            \right.
\end{align}
If $\xi=-1$, we have
\begin{align}\label{eq:ErrAsyBII}
f(x) - f_n(x) &= \mathcal{B}(\alpha,\xi) \left\{
            \begin{array}{ll}
{\displaystyle \frac{1}{2\alpha}n^{-2\alpha} - \frac{1}{2}n^{-2\alpha-1} + O(n^{-2\alpha-2}) },  & \hbox{$x=\xi$,}   \\[12pt]
{\displaystyle -\frac{U_{2n}(\cos(\varphi_{\xi}^{-}))}{2}
n^{-2\alpha-1} + O(n^{-2\alpha-2}) }, & \hbox{$x\neq\xi$.}
            \end{array}
            \right.
\end{align}
\end{itemize}
Here $\mathcal{I}_1(\alpha,\xi)$ and $\mathcal{B}(\alpha,\xi)$ are
defined as in Lemma \ref{lem:AsymChebCoeff}.
\end{theorem}
\begin{proof}
First of all, it is easily seen that $f$ is H\"{o}lder continuous
with exponent $\alpha$ whenever $\alpha\in(0,1)$ and is absolutely
continuous on $\Omega$ whenever $\alpha\geq1$, and hence the
Chebyshev projection $f_n$ converges uniformly to $f$ as
$n\rightarrow\infty$. We now consider the error estimate of $f_n$
and we start with the case $\xi\in(-1,1)$. Due to the uniform
convergence of $f_n$ to $f$, it follows from Lemma
\ref{lem:AsymChebCoeff} that
\begin{align}
f(x) - f_n(x) &= \sum_{k=n+1}^{\infty} a_k T_k(x) \nonumber \\
&= \mathcal{I}_1(\alpha,\xi)\sum_{k=n+1}^{\infty} \frac{T_k(\xi)
T_k(x)}{k^{\alpha+1}} +
\mathcal{I}_2(\alpha,\xi)\sum_{k=n+1}^{\infty} \frac{U_{k-1}(\xi)
T_k(x)}{k^{\alpha+2}} + O(n^{-\alpha-2})  \nonumber \\
&= \frac{\mathcal{I}_1(\alpha,\xi)}{2} \sum_{k=n+1}^{\infty}
\frac{\cos(2k\varphi_{\xi}^{+}) +
\cos(2k\varphi_{\xi}^{-})}{k^{\alpha+1}} \nonumber \\
& \quad  +
\frac{\mathcal{I}_2(\alpha,\xi)}{2\sqrt{1-\xi^2}}\sum_{k=n+1}^{\infty}
\frac{\sin(2k\varphi_{\xi}^{+}) -
\sin(2k\varphi_{\xi}^{-})}{k^{\alpha+2}} + O(n^{-\alpha-2}).
\nonumber
\end{align}
Furthermore, by Lemma \ref{lem:DecaySeries}, the error of $f_n$ can
be rewritten as
\begin{align}\label{eq:LeadTermProj}
f(x) - f_n(x) &= \frac{\mathcal{I}_1(\alpha,\xi)}{2} \left[
\Psi_{\alpha}^{\mathrm{C}}(2\varphi_{\xi}^{+},n) +
\Psi_{\alpha}^{\mathrm{C}}(2\varphi_{\xi}^{-},n) \right]
\nonumber \\[6pt]
& + \frac{\mathcal{I}_2(\alpha,\xi)}{2\sqrt{1-\xi^2}} \left[
\Psi_{\alpha+1}^{\mathrm{S}}(2\varphi_{\xi}^{+},n) -
\Psi_{\alpha+1}^{\mathrm{S}}(2\varphi_{\xi}^{-},n) \right] +
O(n^{-\alpha-2}),
\end{align}
where $\Psi_{\alpha}^{\mathrm{C}}(z,n)$ and
$\Psi_{\alpha}^{\mathrm{S}}(z,n)$ are defined in \eqref{def:PsiFun}.
Let us now consider the pointwise error behavior of Chebyshev
projections. In the case where $x\neq\xi$, it is easy to check that
$\varphi_{\xi}^{+}\in(0,\pi)$ and
$\varphi_{\xi}^{-}\in(-\pi/2,0)\cup(0,\pi/2)$ for all $x\in[-1,1]$.
Applying Lemma \ref{lem:DecaySeries} to \eqref{eq:LeadTermProj}
gives
\begin{align}
f(x)- f_n(x) &= -\frac{\mathcal{I}_1(\alpha,\xi)}{4} \left[
\frac{\sin((2n+1)\varphi_{\xi}^{+})}{\sin(\varphi_{\xi}^{+})} +
\frac{\sin((2n+1)\varphi_{\xi}^{-})}{\sin(\varphi_{\xi}^{-})}
\right] n^{-\alpha-1} + O(n^{-\alpha-2}) \nonumber \\[5pt]
&= -\frac{\mathcal{I}_1(\alpha,\xi)}{4} \left[
U_{2n}(\cos(\varphi_{\xi}^{+})) + U_{2n}(\cos(\varphi_{\xi}^{-}))
\right] n^{-\alpha-1} + O(n^{-\alpha-2}). \nonumber
\end{align}
This proves \eqref{eq:AsyCaseI}. We next consider the case where
$x=\xi$. It is easily seen that $\varphi_{\xi}^{+}=\arccos(\xi)$ and
$\varphi_{\xi}^{-}=0$. Hence, applying Lemma \ref{lem:DecaySeries}
to \eqref{eq:LeadTermProj} again, we obtain that
\begin{align}
f(x)- f_n(x) &= \frac{\mathcal{I}_1(\alpha,\xi)}{2} \big[
\Psi_{\alpha}^{\mathrm{C}}(2\arccos(\xi),n) +
\Psi_{\alpha}^{\mathrm{C}}(0,n) \big] +
O(n^{-\alpha-2})  \nonumber \\[6pt]
&= \frac{\mathcal{I}_1(\alpha,\xi)}{2\alpha} n^{-\alpha} -
\mathcal{I}_1(\alpha,\xi)\frac{U_{2n}(\xi)+1}{4} n^{-\alpha-1} +
O(n^{-\alpha-2}).  \nonumber
\end{align}
This proves \eqref{eq:AsyCaseII}. The pointwise error estimates of
$f_n$ for the case $\xi=\pm1$ can be analyzed in a similar way. In
the case where $\xi=1$, we see that
$\varphi_{\xi}^{+}=\varphi_{\xi}^{-}=\arccos(x)/2\in[0,\pi/2]$.
Applying Lemma \ref{lem:DecaySeries} again we obtain
\begin{align}
f(x) - f_n(x) &= \mathcal{B}(\alpha,\xi) \left\{
            \begin{array}{ll}
{\displaystyle -\frac{U_{2n}(\cos(\varphi_{\xi}^{+}))}{2}
n^{-2\alpha-1} + O(n^{-2\alpha-2})},    &~ \hbox{$x\neq\xi$,}   \\[12pt]
{\displaystyle \frac{1}{2\alpha}n^{-2\alpha} -
\frac{1}{2}n^{-2\alpha-1} + O(n^{-2\alpha-2})},  & \hbox{ $x=\xi$.}
            \end{array}
            \right. \nonumber
\end{align}
This proves the result for $\xi=1$. In the case where $\xi=-1$, we
see that $\varphi_{\xi}^{\pm}=(\arccos(x)\pm\pi)/2$ and
$\varphi_{\xi}^{+}\in[\pi/2,\pi]$ and
$\varphi_{\xi}^{-}\in[-\pi/2,0]$, and therefore
\begin{align}
f(x) - f_n(x) &= \mathcal{B}(\alpha,\xi) \left\{
            \begin{array}{ll}
{\displaystyle
-\frac{U_{2n}(\cos(\varphi_{\xi}^{-}))}{2}n^{-2\alpha-1} +
O(n^{-2\alpha-2})},    & \hbox{$x\neq\xi$,}   \\[12pt]
{\displaystyle \frac{1}{2\alpha} n^{-2\alpha} -
\frac{1}{2}n^{-2\alpha-1} + O(n^{-2\alpha-2})},  & \hbox{$x=\xi$.}
            \end{array}
            \right.  \nonumber
\end{align}
This proves the result for $\xi=-1$ and the proof of Theorem
\ref{thm:PointErrorProj} is complete.
\end{proof}

\begin{figure}[ht]
\centering
\includegraphics[width=0.5\textwidth]{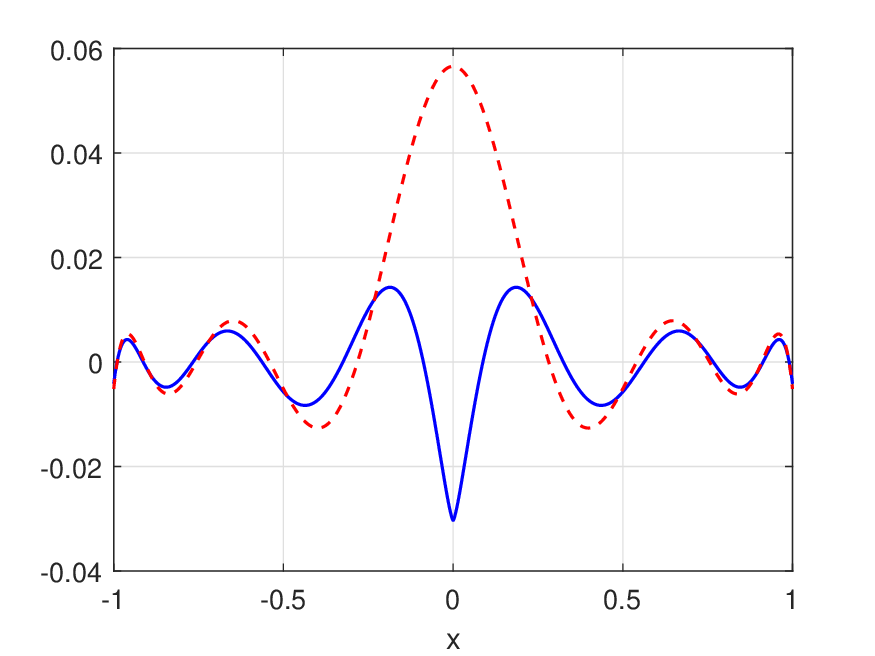}~
\includegraphics[width=0.5\textwidth]{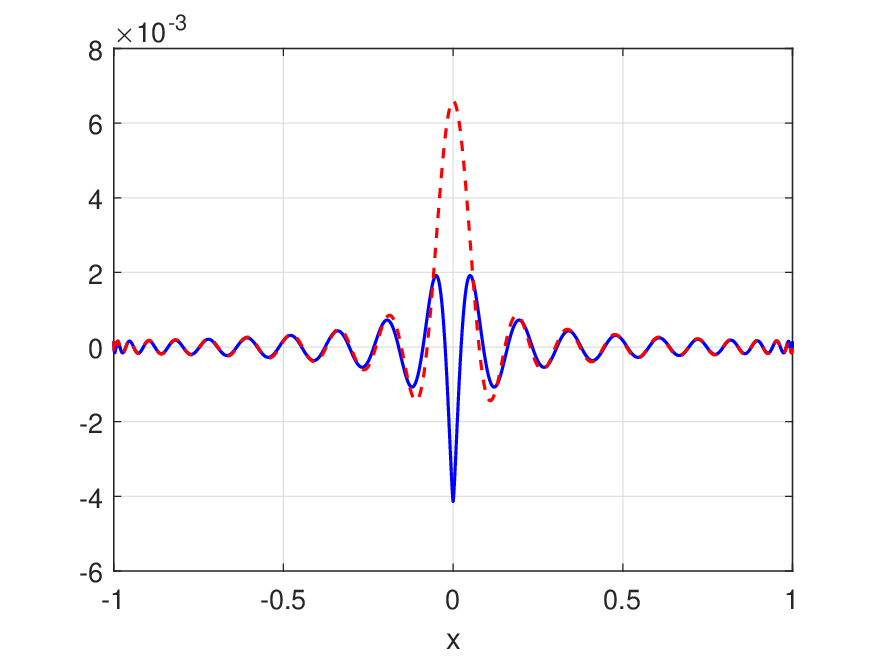}
\caption{The pointwise error curves (line) and the first term of the
asymptotic expansion of the error (dashed) for $n=10$ (left) and
$n=30$ (right). Here $f(x)=|x|^{3/2} e^{\cos(x)}$. }
\label{fig:Envelope}
\end{figure}

To check the pointwise error estimates stated in Theorem
\ref{thm:PointErrorProj}, we plot the pointwise error curves of
$f_n(x)$ and the leading term on the right-hand side of
\eqref{eq:AsyCaseI} in Figure \ref{fig:Envelope} for the function
$f(x)=|x|^{3/2}e^{\cos(x)}$. Clearly, we see that both curves
coincide more closely as $n$ increases except in the neighborhood of
the singularity, which confirms the validity of Theorem
\ref{thm:PointErrorProj}.

\begin{remark}
As $n\rightarrow\infty$, combining
\cite[Theorem~8.21.2]{szego1975orth} with the fact that
$U_n(\pm1)=(\pm1)^n(n+1)$ we obtain
\begin{align}
U_n(x) = \left\{
\begin{array}{ll}
{\displaystyle O(n)},  & \hbox{$x=\pm1$,}   \\[5pt]
{\displaystyle O(1)},  & \hbox{$x\in(-1,1)$.}
\end{array}
\right. \nonumber
\end{align}
Consequently, we can deduce from Theorem \ref{thm:PointErrorProj}
that the rate of convergence of $f_n$ at each $x\neq\xi$ is one
power of $n$ faster than the rate of convergence of $f_n$ at
$x=\xi$, and this justifies the error localization of $f_n$.
\end{remark}

In what follows, we state the error estimate of best approximations
in the maximum norm and its proof can be found, e.g., in
\cite[Chapter~7]{timan1963}.
\begin{theorem}\label{thm:BestRate}
Let $f$ be the function defined in \eqref{def:Model} and let
$p_n^{*}$ be its best approximation of degree $n$. Then, for
$n\gg1$, it holds that
\begin{align}\label{eq:BestRate}
\|f - p_n^{*}\|_{L^{\infty}(\Omega)} = \left\{
\begin{array}{ll}
O(n^{-\alpha}),     & \hbox{if $\xi\in(-1,1)$,}   \\[10pt]
O(n^{-2\alpha}),    & \hbox{if $\xi=\pm1$.}
            \end{array}
            \right.
\end{align}
\end{theorem}
As a consequence, we can conclude from Theorem
\ref{thm:PointErrorProj} and Theorem \ref{thm:BestRate} as follows.
\begin{itemize}
\item In the case of $\xi\in(-1,1)$, the rate of convergence of $p_n^{*}$ is
$O(n^{-\alpha})$ and, owing to the equioscillation property, this
convergence rate is achieved uniformly over the whole interval
$\Omega$. We therefore deduce that $f_n$ converges one power of $n$
faster than $p_n^{*}$ whenever $x$ is away from $\xi$ and $f_n$ and
$p_n^{*}$ converge at the same rate whenever $x=\xi$, and this gives
a clear explanation for those observations displayed in Figure
\ref{fig:Exam}. Moreover, similar conclusions also hold for the case
$\xi=\pm1$.

\item Bernstein in \cite{bernstein1913sur,bernstein1938sur} initiated the
study of the error of the best approximation to the function
$f(x)=|x|^{\alpha}$ and he proved that the limit
\[
\lim_{n\rightarrow\infty} n^{\alpha} \|f -
p_n^{*}\|_{L^{\infty}(\Omega)} := \mu(\alpha)
\]
exists for each $\alpha>0$ and $\alpha\neq2,4,6,\ldots$. However, an
explicit formula for $\mu(\alpha)$ is still unknown. In the special
case $\alpha=1$, it has been shown in \cite{varga1985bern} that
$\mu(\alpha)=0.280169499\ldots$ by using high-precision
calculations. An immediate question that comes to mind is what is
the corresponding result for Chebyshev projections. To address this,
from Theorem \ref{thm:PointErrorProj} we can deduce immediately that
\begin{align}\label{eq:ChebLeadCoeff}
\lim_{n\rightarrow\infty} n^{\alpha} \|f -
f_n\|_{L^{\infty}(\Omega)} = \frac{2\Gamma(\alpha)}{\pi}
\left|\sin\left(\frac{\alpha\pi}{2}\right)\right|.
\end{align}
In the case $\alpha=1$, i.e., $f(x)=|x|$, we see that
$\lim_{n\rightarrow\infty} n^{\alpha} \|f -
f_n\|_{L^{\infty}(\Omega)} = 2/\pi$, and hence the maximum error of
$p_n^{*}$ is better than that of $f_n$ by a factor of about $2.27$
as $n\rightarrow\infty$.
\end{itemize}

\section{Extensions}\label{sec:Extension}
In this section we extend the framework of our analysis to several
closely related topics, including superconvergent points of
Chebyshev projections and pointwise error estimates of Chebyshev
inerpolants, Chebyshev spectral differentiation and Legendre
projections. We will show that the error localization of Chebyshev
inerpolants, Chebyshev spectral differentiation and Legendre
projections can also be justified using similar arguments.

\subsection{Superconvergence points of Chebyshev projections}
Superconvergence theory has received much attention in diverse areas
ranging from $h$-version of finite element methods to collocation
methods for Volterra integral equations (see, e.g.,
\cite{brunner2004coll,superconv1995super}). In the case of spectral
methods, however, only a few studies had been conducted in the
literature. For instance, superconvergence points of spectral
interpolation for functions analytic in a region containing $\Omega$
have been studied in \cite{wang2014super,zhang2012super} and it was
shown that the rate of convergence of Chebyshev and Jacobi-type
interpolants can be improved by some algebraic factors at their
superconvergence points. For functions with limited regularity,
however, the rate of convergence of these spectral interpolants at
their superconvergence points can not be improved anymore.

In the following, we extend our discussion to the superconvergence
points of Chebyshev projections, that is, we seek a set of points
$\{y_{j}^{n}\}\subseteq\Omega$ such that
\begin{align}\label{def:supconv}
n^{\nu} |f(y_{j}^{n})-f_n(y_{j}^{n})| \leq K \|f -
f_n\|_{L^{\infty}(\Omega)},
\end{align}
where $\nu>0$ and $K$ is a generic positive constant independent of
$n$. With \eqref{def:supconv}, one can deduce from Theorem
\ref{thm:PointErrorProj} that $\nu=1$ for all
$x\in\Omega\backslash\{\xi\}$. Below, we will explore the case of
$\nu=2$, that is, we seek a set of points
$\{y_{j}^{n}\}\subseteq\Omega$ such that the rate of convergence of
$f_n$ at these points is faster than the maximum error of $f_n$ by a
factor of $n^2$.

\begin{figure}[ht]
\centering
\includegraphics[width=0.5\textwidth]{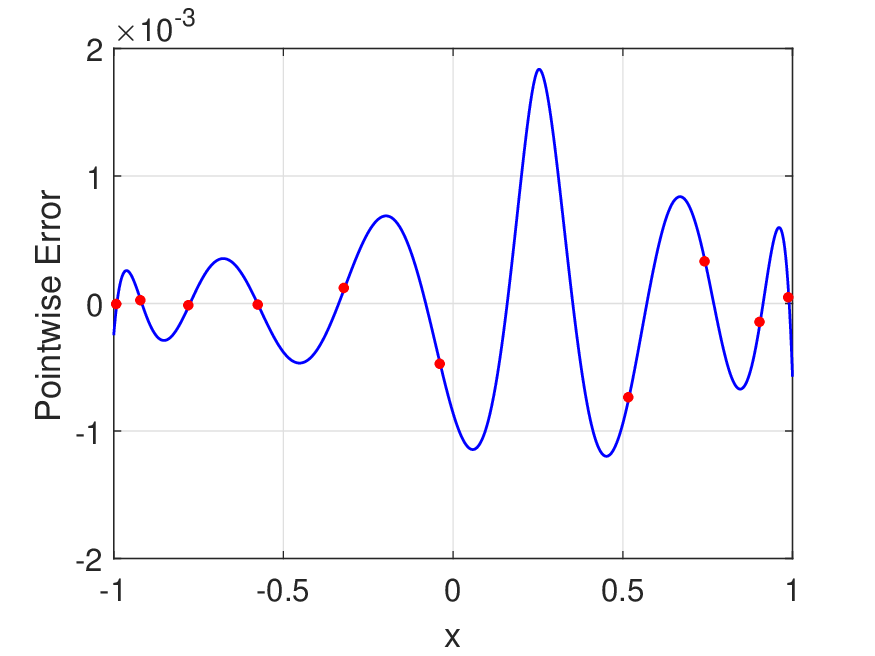}~
\includegraphics[width=0.5\textwidth]{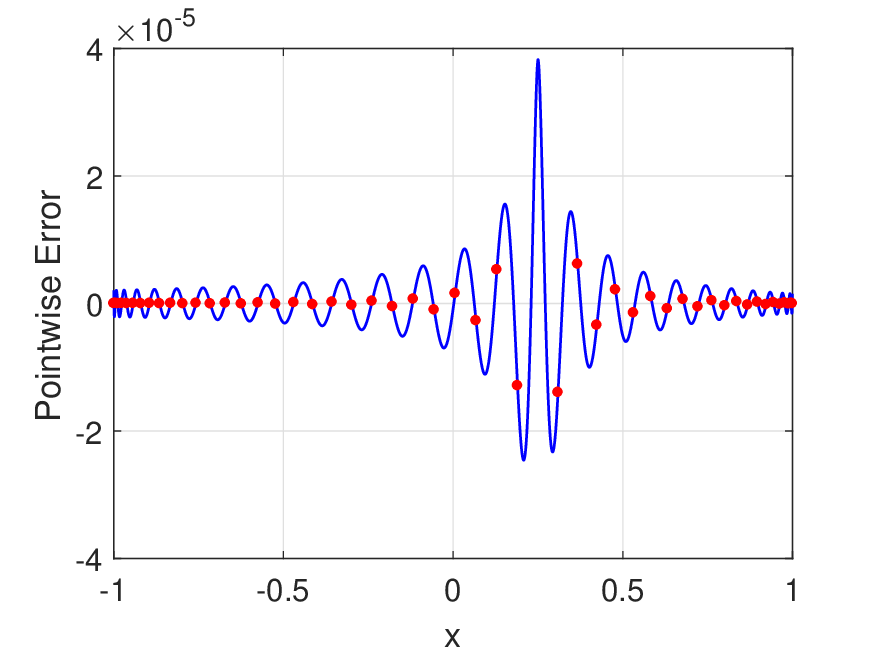}
\caption{The pointwise error curves and the error at the
superconvergent points (dots) for $n=10$ (left) and $n=50$ (right).
Here $f(x)=|x-0.25|^{5/2} e^{x}$. } \label{fig:ExamS}
\end{figure}

\begin{figure}[ht]
\centering
\includegraphics[width=0.5\textwidth]{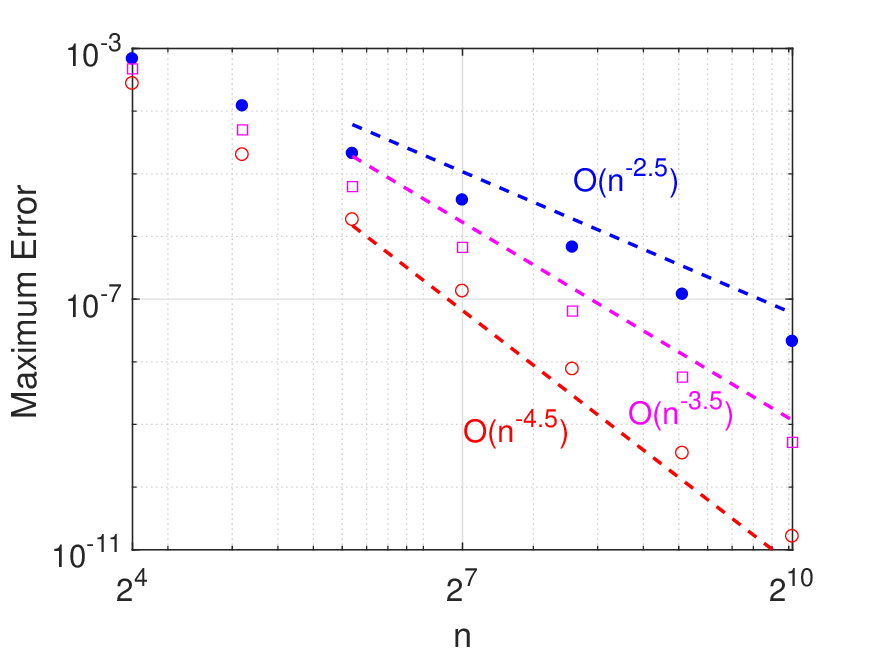}~
\includegraphics[width=0.5\textwidth]{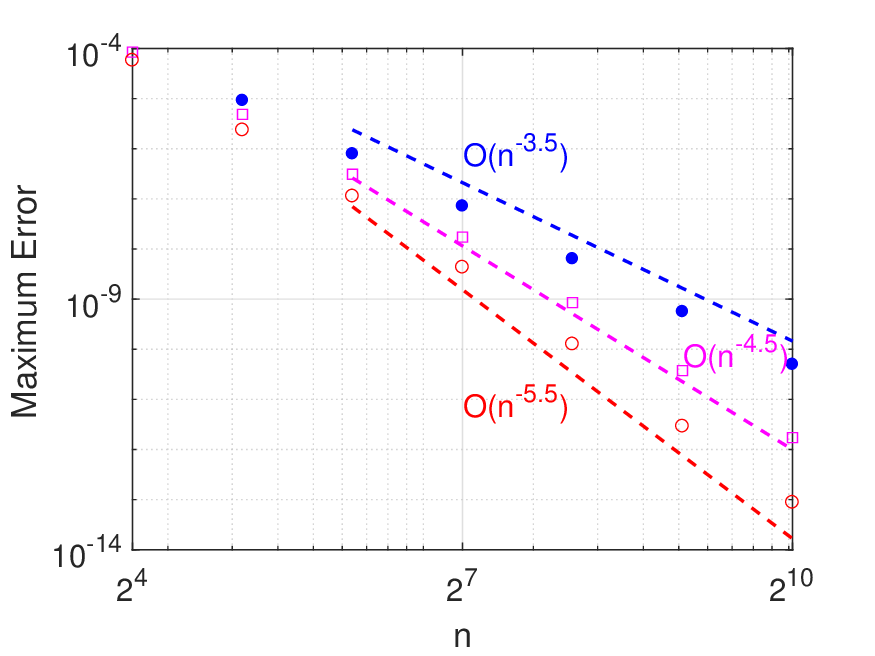}
\caption{Maximum errors of $f_n$ on $\Omega$ (dots) and on
$\Omega_{\varepsilon}(\xi)$ (box) and at the set
$\mathcal{S}_{\varepsilon}(\xi)$ (circles) for $\alpha=2.5$ (left)
and $\alpha=3.5$ (right). Here $\xi=0.25$, $g(x)=e^{x}$ and
$\varepsilon=0.1$. } \label{fig:ExamI}
\end{figure}

We start with the case $\xi\in(-1,1)$. Firstly, from Theorem
\ref{thm:PointErrorProj} we know that the maximum error of $f_n$ on
$\Omega$ satisfies $\|f - f_n
\|_{L^{\infty}(\Omega)}=O(n^{-\alpha})$ and this error estimate is
attained at the singularity $x=\xi$. To obtain superconvergence
points which satisfy \eqref{def:supconv} with $\nu=2$, we impose the
following condition
\begin{align}\label{eq:SuperPoint}
U_{2n}(\cos(\varphi_{\xi}^{+})) + U_{2n}(\cos(\varphi_{\xi}^{-})) =
0,
\end{align}
where, using elementary algebraic manipulations,
\begin{align}\label{eq:cosvar}
\cos(\varphi_{\xi}^{\pm}) =
\frac{\sqrt{(1+x)(1+\xi)}\mp\sqrt{(1-x)(1-\xi)}}{2}.
\end{align}
It is easily verified that the left-hand side of
\eqref{eq:SuperPoint} is a polynomial of degree at most $n$, from
which we can obtain a set of superconvergence points
$\{y_{j}^{n}(\xi)\}$. In Figure \ref{fig:ExamS} we illustrate the
error curves of $f_n$ and the errors at the superconvergence points.
We can see that the error of $f_n$ at the superconvergence points
are smaller than the maximum error of $f_n$, especially when the
superconvergence points are away from the singularity. Moreover, we
also observe that the larger $n$, the more superconvergence points
distributed in the neighborhood of the singularity $x=\xi$. Recall
that the maximum error of $f_n$ is always attained at the
singularity for moderate and large values of $n$. Hence, to
guarantee \eqref{def:supconv}, we introduce a small parameter
$\varepsilon$ which is independent of $n$ and $0<\varepsilon\ll1$.
We denote by $\Omega_{\varepsilon}(\xi)$ the interval $\Omega$
except the $\varepsilon$-neighborhood of $\xi$, i.e.,
$\Omega_{\varepsilon}(\xi)=\{x\in\Omega:|x-\xi|\geq\varepsilon\}$
and denote by $S(\xi)$ the solution set of \eqref{eq:SuperPoint}. We
then define the set of superconvergence points over the interval
$\Omega_{\varepsilon}(\xi)$, i.e.,
\begin{align}
S_{\varepsilon}(\xi) = \left\{ y_{j}^{n}(\xi)\in
S(\xi)\cap\Omega_{\varepsilon}(\xi) \right\}.
\end{align}
As a consequence of Theorem \ref{thm:PointErrorProj}, the rate of
convergence of $f_n$ at the points in $S_{\varepsilon}(\xi)$ is
$O(n^{-\alpha-2})$, which is faster than the maximum error of $f_n$
over the interval $\Omega$ by a factor of $n^2$ and is faster than
the maximum error of $f_n$ over the interval
$\Omega_{\varepsilon}(\xi)$ by a factor of $n$. In Figure
\ref{fig:ExamI} we illustrate the maximum errors of $f_n$ on
$\Omega$ and on $\Omega_{\varepsilon}(\xi)$ and at the set
$S_{\varepsilon}(\xi)$ for two values of $\alpha$. Clearly, we see
that numerical results are in good agreement with our theoretical
analysis.

\begin{remark}
In the special case $\xi=0$, the equation \eqref{eq:SuperPoint} can
be solved explicitly. Specifically, from \eqref{eq:cosvar} we know
that $\cos(\varphi_{\xi}^{\pm})=(\sqrt{1+x}\mp\sqrt{1-x})/2$ and
thus $\cos(2\varphi_{\xi}^{\pm}) = \mp\sqrt{1-x^2}$. Let $W_n(x)$ be
the Chebyshev polynomial of the fourth kind of degree $n$ (see,
e.g., \cite[Chapter~1]{mason2003cheb}). From the definitions of
$U_n(x)$ and $W_n(x)$, we obtain that
\begin{align}\label{eq:SuperZeroI}
U_{2n}(\cos(\varphi_{\xi}^{+})) + U_{2n}(\cos(\varphi_{\xi}^{-})) &=
W_{n}(\cos(2\varphi_{\xi}^{+})) + W_{n}(\cos(2\varphi_{\xi}^{-})) \nonumber \\
&= W_{n}(-\sqrt{1-x^2}) + W_{n}(\sqrt{1-x^2}).
\end{align}
Moreover, from \cite[Equation~(1.18)]{mason2003cheb} we know that
$W_n(x)+W_n(-x)=2U_n(x)$ whenever $n$ is even and
$W_n(x)+W_n(-x)=2U_{n-1}(x)$ whenever $n$ is odd. Combining this
with \eqref{eq:SuperZeroI} we can deduce that the solution of
\eqref{eq:SuperPoint} is
$S(\xi)=\{y_{j}^{n}(\xi)=\sin(j\pi/(n+1)),j=1,\ldots,n \}$ whenever
$n$ is even and $S(\xi)=\{y_{j}^{n}(\xi)=\sin(j\pi/n),j=1,\ldots,n-1
\}$ whenever $n$ is odd.
\end{remark}

\begin{figure}[ht]
\centering
\includegraphics[width=0.5\textwidth]{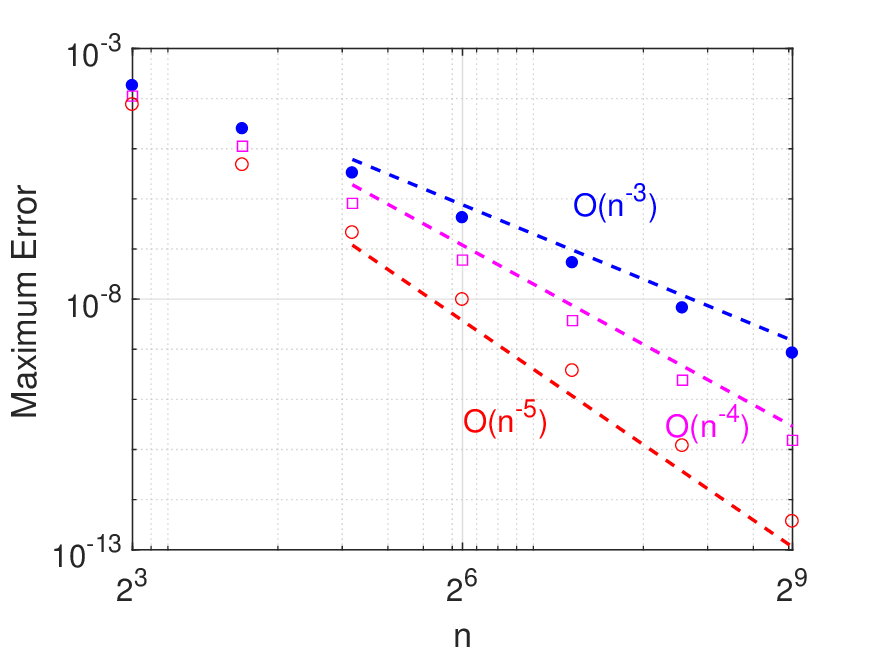}~
\includegraphics[width=0.5\textwidth]{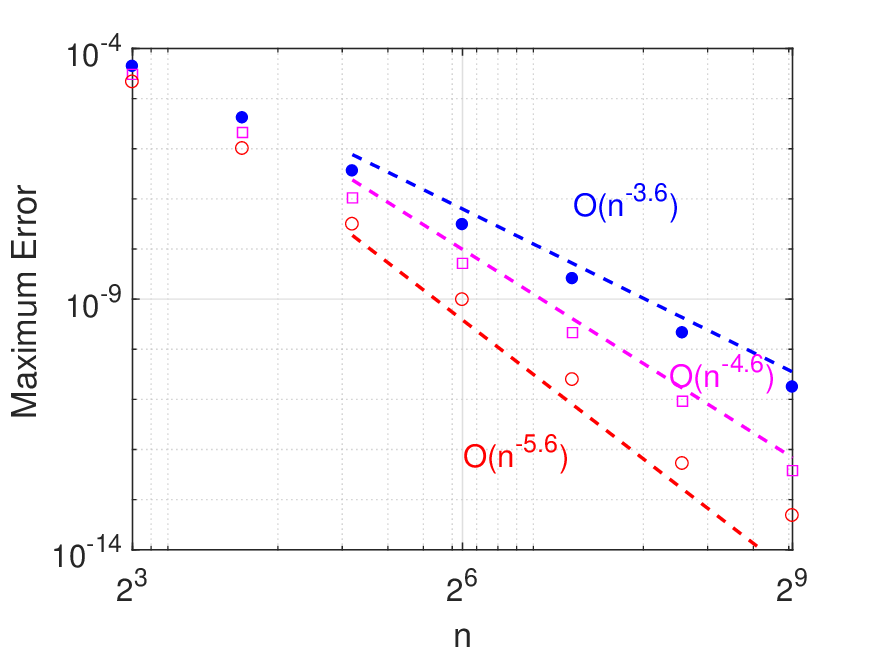}
\caption{Maximum errors of $f_n$ (dots) on $\Omega$ and on
$\Omega_{\varepsilon}(\xi)$ (box) and at the superconvergence points
$S_{\varepsilon}(\xi)$ (circles) for $\alpha=1.5$ (left) and
$\alpha=1.8$ (right). Here $\xi=-1$, $g(x)=1/(3-x)$ and
$\varepsilon=0.05$.} \label{fig:ExamII}
\end{figure}

Now we turn to the case $\xi=\pm1$. It is known that the maximum
error of $f_n$ on $\Omega$ satisfies $\|f - f_n
\|_{L^{\infty}(\Omega)}=O(n^{-2\alpha})$ and this error estimate is
attained at $x=\xi$. In this case, the superconvergence points can
be obtained by imposing the condition
$U_{2n}(\cos(\varphi_{\xi}^{-}))=0$ whenever $\xi=-1$ and
$U_{2n}(\cos(\varphi_{\xi}^{+}))=0$ whenever $\xi=1$, from which we
obtain
\begin{align}\label{eq:supend}
S(\xi) = \left\{ y_{j}^{n}(\xi) = \xi\cos\left(\frac{2\pi j}{2n+1}
\right), ~~j = 1,\ldots,n \right\}.
\end{align}
Moreover, it is easily verified that these points are zeros of
$V_n(x)$ whenever $\xi=-1$ and of $W_n(x)$ whenever $\xi=1$, where
$V_n(x)$ and $W_n(x)$ denote respectively the Chebyshev polynomials
of the third and fourth kind of degree $n$. Let
$\Omega_{\varepsilon}(\xi)$ and $S_{\varepsilon}(\xi)$ be defined as
before. It is easily verified that
\[
|y_{j}^{n}(\xi) - \xi| \geq \varepsilon ~~~ \Longrightarrow ~~~
j\geq \left\lceil \frac{2n+1}{2\pi} \arccos(1-\varepsilon)
\right\rceil :=n_{\varepsilon},
\]
where $\lceil x\rceil$ denotes the ceiling of $x$, then we obtain
$S_{\varepsilon}(\xi)=\{y_{j}^{n}(\xi),j=n_{\varepsilon},\ldots,n
\}$ where $y_{j}^{n}(\xi)$ is defined as in \eqref{eq:supend}. From
Theorem \ref{thm:PointErrorProj} we can deduce that the rate of
convergence of $f_n$ at the superconvergence points in
$S_{\varepsilon}(\xi)$ is $O(n^{-2\alpha-2})$, which is faster than
the maximum error of $f_n$ over the interval $\Omega$ by a factor of
$n^2$ and is faster than the maximum error of $f_n$ over the
interval $\Omega_{\varepsilon}(\xi)$ by a factor of $n$. In Figure
\ref{fig:ExamII} we illustrate the maximum errors of $f_n$ on
$\Omega$ and on $\Omega_{\varepsilon}(\xi)$ and at these
superconvergence points $S_{\varepsilon}(\xi)$ for two different
values of $\alpha$. Clearly, we see that these numerical results are
in accordance with our theoretical analysis.

\subsection{Chebyshev interpolants}
In this subsection we consider pointwise error estimates of
Chebyshev interpolants, i.e., polynomial interpolants at the zeros
or extrema of Chebyshev polynomials. As will be shown below, the
results are analogous to those of Theorem \ref{thm:PointErrorProj},
but their derivation is a little more involved.

\subsubsection{Chebyshev points of the first kind}
Let $\{x_j\}_{j=0}^{n}$ be the Chebyshev points of the first kind,
i.e., $x_j = \cos((j+1/2)\pi/(n+1))$, and let $p_n^{\mathrm{I}}$ be
the unique polynomial of degree $n$ which interpolates $f$ at
$\{x_j\}_{j=0}^{n}$. For $k=0,\ldots,n$, by straightforward
calculations, one can verify that
\begin{align}\label{def:DisOrthI}
\sum_{s=0}^{n} T_j(x_s) T_k(x_s) = (-1)^{\ell} \left\{
            \begin{array}{ll}
{\displaystyle n+1},           &~~ \hbox{$k=0,~~j=2\ell(n+1)$, $\ell\in\mathbb{N}_0$,} \\[8pt]
{\displaystyle \frac{n+1}{2}}, &~~ \hbox{$k\in \mathcal{U}_n$, $j=2\ell(n+1)+k$, $\ell\in\mathbb{N}_0$, } \\[8pt]
{\displaystyle \frac{n+1}{2}}, &~~ \hbox{$k\in \mathcal{U}_n$, $j=2\ell(n+1)-k$, $\ell\in\mathbb{N}$, } \\[8pt]
{\displaystyle 0}, &~~ \hbox{otherwise,}
            \end{array}
            \right.
\end{align}
where $\mathcal{U}_n=\{1,\ldots,n\}$. Based on the discrete
orthogonality \eqref{def:DisOrthI}, the Chebyshev interpolant
$p_n^{\mathrm{I}}(x)$ can be written explicitly as
\begin{align}\label{def:ChebInterpI}
p_n^{\mathrm{I}}(x) = \sum_{k=0}^{n}{'} b_k T_k(x), \quad  b_k =
\frac{2}{n+1} \sum_{j=0}^{n} f(x_j) T_k(x_j),
\end{align}
and it is known that these coefficients $\{b_0,\ldots,b_n\}$ can be
computed rapidly by using the FFT. We now turn to the task of
analyzing the pointwise error estimate of $p_n^{\mathrm{I}}(x)$ and
the key ingredient of our analysis is the aliasing formula of
Chebyshev coefficients. Specifically, plugging the Chebyshev series
\eqref{def:ChebExp} into $b_k$ and applying \eqref{def:DisOrthI}, we
obtain
\begin{align}
b_k = \frac{2}{n+1} \sum_{j=0}^{n} \sum_{\ell=0}^{\infty}{'}
a_{\ell} T_{\ell}(x_j) T_k(x_j) & = \frac{2}{n+1}
\sum_{\ell=0}^{\infty}{'} a_{\ell} \sum_{j=0}^{n}
T_{\ell}(x_j) T_k(x_j) \nonumber \\
&= a_k + \sum_{\ell=1}^{\infty} (-1)^{\ell} \left( a_{2\ell(n+1)-k}
+ a_{2\ell(n+1)+k} \right). \nonumber
\end{align}
Combining this with \eqref{def:ChebExp} and \eqref{def:ChebInterpI},
we get
\begin{align}\label{eq:ChebErrorS1}
f(x) - p_n^{\mathrm{I}}(x) &= \frac{a_0-b_0}{2} + \sum_{k=1}^{n}
(a_k - b_k) T_k(x) + \sum_{k=n+1}^{\infty} a_k T_k(x) \nonumber \\
&= \sum_{k=n+1}^{\infty} a_k \big( T_k(x) + \nu(k) T_{\eta(k)}(x)
\big),
\end{align}
where $\eta(k)=|(k+n)(\mathrm{mod}~ 2(n+1)) - n|$ and
\begin{align}
\nu(k)= \left\{
            \begin{array}{ll}
{\displaystyle   0},  &~~ \hbox{$k=(2\ell+1)(n+1)$, $\ell\in\mathbb{N}_0$,} \\[8pt]
{\displaystyle   (-1)^{\lfloor (k-n-1)/(2n+2) \rfloor}},  &~~
\hbox{otherwise.}
            \end{array}
            \right. \nonumber
\end{align}
It is easily seen that $\{\eta(k)\}$ is an $(2n+2)$-periodic
sequence, i.e., $\eta(k)=\eta(k+2n+2)$ for each $k\geq n+2$. In the
following, for simplicity of notation, we set
$\theta=\arccos(x)\in[0,\pi]$. By using the periodic property of
$\{\eta(k)\}$, and after somewhat tedious but quite elementary
calculations, we arrive at
\begin{align}\label{eq:ChebInterpErrorI}
f(x) - p_n^{\mathrm{I}}(x) &= \sum_{k=n+1}^{\infty} a_k \cos(k\theta) \nonumber \\
&~ + \sum_{\ell=1}^{\infty} \cos(2(2\ell-1)(n+1)\theta)
\sum_{k=2(2\ell-1)(n+1)-n}^{2(2\ell-1)(n+1)+n} a_k\cos(k\theta)
\nonumber \\
&~ + \sum_{\ell=1}^{\infty} \sin(2(2\ell-1)(n+1)\theta)
\sum_{k=2(2\ell-1)(n+1)-n}^{2(2\ell-1)(n+1)+n} a_k\sin(k\theta)
\nonumber \\
&~ - \sum_{\ell=1}^{\infty} \cos(4\ell(n+1)\theta)
\sum_{k=4\ell(n+1)-n}^{4\ell(n+1)+n} a_k\cos(k\theta)  \\
&~ - \sum_{\ell=1}^{\infty} \sin(4\ell(n+1)\theta)
\sum_{k=4\ell(n+1)-n}^{4\ell(n+1)+n} a_k\sin(k\theta). \nonumber
\end{align}
Notice that the first sum on the right-hand side of
\eqref{eq:ChebInterpErrorI} is exactly the remainder term of $f_n$,
its estimate follows immediately from Theorem
\ref{thm:PointErrorProj}. For the last four sums on the right-hand
side of \eqref{eq:ChebInterpErrorI}, their asymptotic estimates
follow by means of Lemma \ref{lem:AsymChebCoeff} and Lemma
\ref{lem:DecaySeries}. We omit the details of the derivation here
since it is lengthy but similar to that of Theorem
\ref{thm:PointErrorProj}, and give the final results below: Let
$\mathcal{N}(\xi)$ be a small neighborhood of the singularity $\xi$
in $\Omega$ and $|\mathcal{N}(\xi)|\rightarrow0$ as
$n\rightarrow\infty$.
\begin{itemize}
\item[\rm (i)] If $\xi\in(-1,1)$ is not an interpolation point, then for $x\neq x_j$ we
have
\begin{align}\label{eq:ChebRateI}
|f(x) - p_n^{\mathrm{I}}(x)| = \left\{
            \begin{array}{ll}
{\displaystyle O(n^{-\alpha}) },    & \hbox{if~$x\in\mathcal{N}(\xi)$,}   \\[8pt]
{\displaystyle O(n^{-\alpha-1}) },  & \hbox{otherwise.}
            \end{array}
            \right.
\end{align}
If $\xi\in(-1,1)$ happens to be an interpolation point, then the
error of $p_n^{\mathrm{I}}(x)$ is zero whenever $x=\xi$ and behaves
like $O(n^{-\alpha})$ whenever $x$ is very close to $\xi$ and
behaves like $O(n^{-\alpha-1})$ whenever $x$ is away from $\xi$.

\item[\rm (ii)] If $\xi=\pm1$, then for $x\neq
x_j$ we have
\begin{align}\label{eq:ChebRateII}
|f(x) - p_n^{\mathrm{I}}(x)| = \left\{
            \begin{array}{ll}
{\displaystyle O(n^{-2\alpha}) },    & \hbox{if~$x\in\mathcal{N}(\xi)$,}   \\[8pt]
{\displaystyle O(n^{-2\alpha-1}) },  & \hbox{otherwise.}
            \end{array}
            \right.
\end{align}
\end{itemize}
Here we remark that the leading coefficients of the error estimates
of $p_n^{\mathrm{I}}(x)$ involve rather lengthy expression and
therefore are omitted. From the above error estimates, we can see
clearly that the pointwise error estimates of $p_n^{\mathrm{I}}(x)$
are quite similar to that of $f_n$ and these results justify the
error localization of $p_n^{\mathrm{I}}(x)$; see Figure
\ref{fig:ExamIII} for an illustration.

\begin{remark}
The accuracy of $p_n^{\mathrm{I}}(x)$ and $p_n^{*}(x)$ was compared
in \cite{li2004optimal} from the point of view of measuring their
maximum errors and it was shown that $p_n^{\mathrm{I}}(x)$ is just
as good as $p_n^{*}(x)$ in computing elementary functions, such as
$f(x)=e^x,\sin(x)$. Here, our results \eqref{eq:ChebRateI} and
\eqref{eq:ChebRateII} provide a more thorough insight into the
comparison of $p_n^{\mathrm{I}}(x)$ and $p_n^{*}(x)$, that is,
$p_n^{\mathrm{I}}(x)$ actually converges faster than $p_n^{*}(x)$ by
one power of $n$ except in a neighborhood of the singularity.
\end{remark}

\subsubsection{Chebyshev points of the second kind}
\begin{figure}[ht]
\centering
\includegraphics[width=0.5\textwidth]{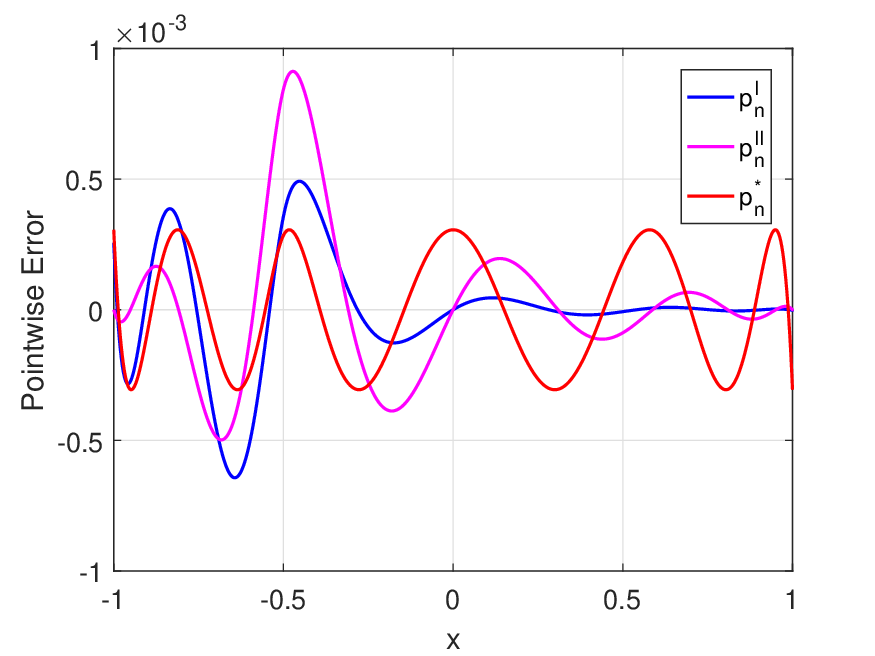}~
\includegraphics[width=0.5\textwidth]{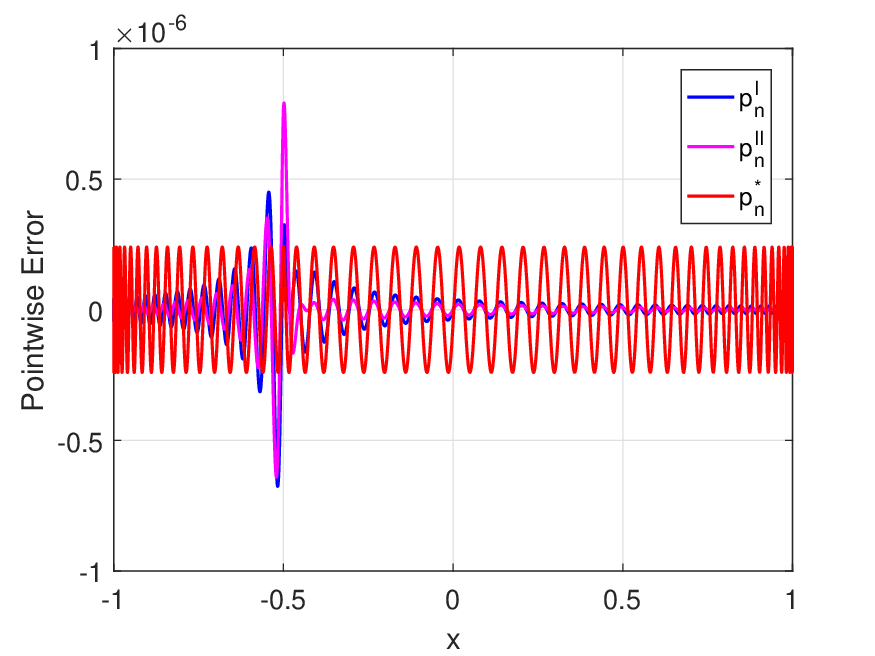}
\caption{Pointwise error curves of $p_n^{\mathrm{I}}(x)$,
$p_n^{\mathrm{II}}(x)$ and $p_n^{*}(x)$ for $f(x)=|x+0.5|^3e^x$ with
$n=10$ (left) and $n=100$ (right). } \label{fig:ExamIII}
\end{figure}

Let $\{x_j\}_{j=0}^{n}$ be the set of Chebyshev points of the second
kind (also known as Chebyshev--Lobatto points or Clenshaw--Curtis
points), i.e., $ x_j=\cos(j\pi/n)$, and let $p_n^{\mathrm{II}}(x)$
be the polynomial of degree $n$ which interpolates $f(x)$ at the
Chebyshev points of the second kind. It is straightforward to verify
that
\begin{align}\label{eq:DisOrthII}
\sum_{s=0}^{n}{''} T_j(x_s) T_k(x_s) = \left\{
            \begin{array}{ll}
{\displaystyle   n},  &~~ \hbox{$(j+k)\in \mathcal{V}_n,~~|j-k|\in \mathcal{V}_n$,} \\[8pt]
{\displaystyle n/2},  &~~ \hbox{$(j+k)\not\in \mathcal{V}_n,~~|j-k|\in \mathcal{V}_n$,} \\[8pt]
{\displaystyle n/2},  &~~ \hbox{$(j+k)\in \mathcal{V}_n,~~|j-k|\not\in \mathcal{V}_n$,} \\[8pt]
{\displaystyle   0},  &~~ \hbox{$(j+k)\not\in \mathcal{V}_n,
~~|j-k|\not\in \mathcal{V}_n$,}
            \end{array}
            \right.
\end{align}
where the double prime indicates that both the first and last terms
of the summation are to be halved and
$\mathcal{V}_n=\{k:k\in\mathbb{N}_{0}~\mbox{and}~
k(\mathrm{mod}~2n)=0\}$, the Chebyshev interpolant
$p_n^{\mathrm{II}}(x)$ can be written explicitly as
\begin{align}\label{eq:ChebInterpII}
p_n^{\mathrm{II}}(x) = \sum_{k=0}^{n}{''} c_k T_k(x), \quad  c_k =
\frac{2}{n} \sum_{j=0}^{n}{''} f(x_j) T_k(x_j).
\end{align}
Moreover, it is well known that these coefficients
$\{c_0,\ldots,c_n\}$ can be computed rapidly by using the FFT in
only $O(n\log n)$ operations. Now we turn to the error estimate of
$p_n^{\mathrm{II}}(x)$ and the derivation is similar to that of
$p_n^{\mathrm{I}}(x)$. Plugging the Chebyshev series into $c_k$ and
applying the discrete orthogonality \eqref{eq:DisOrthII}, we obtain
\begin{align}
c_k = \frac{2}{n} \sum_{j=0}^{n}{''} \sum_{\ell=0}^{\infty}{'}
a_{\ell} T_{\ell}(x_j) T_k(x_j) & = \frac{2}{n}
\sum_{\ell=0}^{\infty}{'} a_{\ell} \sum_{j=0}^{n}{''}
T_{\ell}(x_j) T_k(x_j) \nonumber \\
&= a_k + \sum_{\ell=1}^{\infty}( a_{2\ell n - k} + a_{2\ell n + k}).
\nonumber
\end{align}
The combination of the last equality with \eqref{eq:ChebInterpII}
gives
\begin{align}
f(x) - p_n^{\mathrm{II}}(x) &= \sum_{k=0}^{\infty}{'} a_k T_k(x) -
\sum_{k=0}^{n}{''} c_k T_k(x) = \sum_{k=n+1}^{\infty} a_k \left(
T_k(x) - T_{\psi(k)}(x) \right),
\end{align}
where $\psi(k)=|(k+n-1)(\mathrm{mod}~ 2n)-(n-1)|$. Moreover, it is
easy to verify that $\{\psi(k)\}$ is an $2n$-periodic sequence,
i.e., $\psi(k)=\psi(k+2n)$ for each $k\geq n+1$. By using the
periodic property of $\{\psi(k)\}$, and after some calculations, we
arrive at
\begin{align}\label{eq:ChebInterpErrorII}
f(x) - p_n^{\mathrm{II}}(x) = \sum_{k=n+1}^{\infty} a_k
\cos(k\theta) &- \sum_{\ell=1}^{\infty} \cos(2\ell n\theta)
\sum_{k=(2\ell-1)n+1}^{(2\ell+1)n} a_k \cos(k\theta) \nonumber \\
& - \sum_{\ell=1}^{\infty} \sin(2\ell n\theta)
\sum_{k=(2\ell-1)n+1}^{(2\ell+1)n} a_k \sin(k\theta),
\end{align}
where $\theta=\arccos(x)$. Notice that the first sum is exactly the
remainder term of $f_n$, its estimate follows immediately from
Theorem \ref{thm:PointErrorProj}. For the last two sums on the
right-hand side of \eqref{eq:ChebInterpErrorII}, their asymptotic
estimates can be obtained using Lemma \ref{lem:AsymChebCoeff} and
Lemma \ref{lem:DecaySeries}. Here we omit the details of the
derivation and give the final results below: Let $\mathcal{N}(\xi)$
be defined as in the above subsection.
\begin{itemize}
\item[\rm (i)] If $\xi\in(-1,1)$ is not an interpolation point, then we have for $x\neq
x_j$ that
\begin{align}\label{eq:ChebRateIII}
|f(x) - p_n^{\mathrm{II}}(x)| = \left\{
            \begin{array}{ll}
{\displaystyle O(n^{-\alpha}) },    & \hbox{if~ $x\in\mathcal{N}(\xi)$,}   \\[8pt]
{\displaystyle O(n^{-\alpha-1}) },  & \hbox{otherwise.}
            \end{array}
            \right.
\end{align}
If $\xi\in(-1,1)$ happens to be an interpolation point, then the
error of $p_n^{\mathrm{I}}(x)$ is zero whenever $x=\xi$ and behaves
like $O(n^{-\alpha})$ whenever $x$ is very close to $\xi$ and
behaves like $O(n^{-\alpha-1})$ whenever $x$ is away from $\xi$.

\item[\rm (ii)] If $\xi=\pm1$, then for $x\neq
x_j$ we have
\begin{align}\label{eq:ChebRateIV}
|f(x) - p_n^{\mathrm{II}}(x)| = \left\{
            \begin{array}{ll}
{\displaystyle O(n^{-2\alpha})},     & \hbox{if~$x\in\mathcal{N}(\xi)$ and $x\neq\xi$,}   \\[8pt]
{\displaystyle O(n^{-2\alpha-1}) },  & \hbox{otherwise.}
            \end{array}
            \right.
\end{align}
\end{itemize}
Similar to the case of $p_n^{\mathrm{I}}(x)$, the leading
coefficients of the error estimates of $p_n^{\mathrm{II}}(x)$ also
involve rather lengthy expression and therefore are omitted. We see
that the pointwise error estimates of $p_n^{\mathrm{II}}(x)$ are
quite similar to that of $f_n$ and $p_n^{\mathrm{I}}(x)$; see Figure
\ref{fig:ExamIII}. As a final remark, we point out that our results
\eqref{eq:ChebRateIII} and \eqref{eq:ChebRateIV} provide a
justification for the error localization of $p_n^{\mathrm{II}}(x)$,
which was observed by Trefethen in \cite[Myth~3]{trefethen2011six}.

\subsection{Chebyshev spectral differentiation}
In this subsection, we extend our discussion to the pointwise error
estimate of Chebyshev spectral differentiation. The main result is
stated as follows.
\begin{theorem}\label{thm:ChebDiff}
Let $f$ be the function defined in \eqref{def:Model} for some
$\alpha>1$ and let $f_n$ be its Chebyshev projection of degree $n$.
As $n\rightarrow\infty$, the following pointwise error estimates are
true.
\begin{itemize}
\item[\rm (i)] If $\xi\in(-1,1)$, we have
\begin{align}\label{eq:ChebDiffI}
|f{'}(x) - f_n{'}(x)| &= \left\{
            \begin{array}{ll}
{\displaystyle O(n^{-\alpha}) },    & \hbox{$x\in(-1,1)$,}   \\[8pt]
{\displaystyle O(n^{1-\alpha}) },   & \hbox{$x=\pm1$.}
            \end{array}
            \right.
\end{align}

\item[\rm (ii)] If $\xi=\pm1$, we have
\begin{align}\label{eq:ChebDiffII}
|f{'}(x) - f_n{'}(x)| &= \left\{
            \begin{array}{ll}
{\displaystyle O(n^{-2\alpha}) },   & \hbox{$x\in(-1,1)$,}   \\[8pt]
{\displaystyle O(n^{1-2\alpha}) },  & \hbox{$x=-\xi$,}       \\[8pt]
{\displaystyle O(n^{2-2\alpha}) },  & \hbox{$x=\xi$.}
            \end{array}
            \right.
\end{align}
\end{itemize}
\end{theorem}
\begin{proof}
We only sketch the proof in the case where $\xi\in(-1,1)$. Recalling
the fact that $T_k{'}(x)=kU_{k-1}(x)$ for $k\in\mathbb{N}$ and using
Lemma \ref{lem:AsymChebCoeff}, we obtain for each $x\in(-1,1)$ that
\begin{align}
f{'}(x) - f_n{'}(x) &= \sum_{k=n+1}^{\infty} a_k k U_{k-1}(x) \nonumber \\
&= \mathcal{I}_1(\alpha,\xi) \sum_{k=n+1}^{\infty}
\frac{T_k(\xi)U_{k-1}(x)}{k^{\alpha}} + \mathcal{I}_2(\alpha,\xi)
\sum_{k=n+1}^{\infty}
\frac{U_{k-1}(\xi)U_{k-1}(x)}{k^{\alpha+1}} \nonumber \\
&~~~ + O(n^{-\alpha-1}) \nonumber \\
&= \frac{\mathcal{I}_1(\alpha,\xi)}{2\sqrt{1-x^2}}
\sum_{k=n+1}^{\infty} \frac{\sin(2k\varphi_{\xi}^{+}) +
\sin(2k\varphi_{\xi}^{-})}{k^{\alpha}} \nonumber \\
&~~~ + \frac{\mathcal{I}_2(\alpha,\xi)}{2\sqrt{(1-x^2)(1-\xi^2)}}
\sum_{k=n+1}^{\infty} \frac{\cos(2k\varphi_{\xi}^{-}) -
\cos(2k\varphi_{\xi}^{+})}{k^{\alpha+1}} + O(n^{-\alpha-1})
\nonumber \\
&= \frac{\mathcal{I}_1(\alpha,\xi)}{2\sqrt{1-x^2}}
\left[\Psi_{\alpha-1}^{\mathrm{S}}(2\varphi_{\xi}^{+},n) +
\Psi_{\alpha-1}^{\mathrm{S}}(2\varphi_{\xi}^{-},n)
\right] \nonumber \\
&~~~ + \frac{\mathcal{I}_2(\alpha,\xi)}{2\sqrt{(1-x^2)(1-\xi^2)}}
\left[\Psi_{\alpha}^{\mathrm{C}}(2\varphi_{\xi}^{-},n) -
\Psi_{\alpha}^{\mathrm{C}}(2\varphi_{\xi}^{+},n) \right] +
O(n^{-\alpha-1}), \nonumber
\end{align}
where $\varphi_{\xi}^{+}$ and $\varphi_{\xi}^{-}$ are defined as in
\eqref{def:varphi}. Applying Lemma \ref{lem:DecaySeries} to the
leading two terms on the right-hand side we immediately deduce that
$f{'}(x) - f_n{'}(x)=O(n^{-\alpha})$ for each $x\in(-1,1)$. For
$x=\pm1$, recalling the fact that $U_{k}(\pm1)=(\pm1)^{k}(k+1)$ for
all $k\in\mathbb{N}_0$ we obtain
\begin{align}
f{'}(\pm1) - f_n{'}(\pm1) &= \mathcal{I}_1(\alpha,\xi)
\sum_{k=n+1}^{\infty} \frac{(\pm1)^{k-1}T_k(\xi)}{k^{\alpha-1}} +
O(n^{-\alpha}) = O(n^{1-\alpha}),  \nonumber
\end{align}
where we have used Lemma \ref{lem:DecaySeries} in the last step.
This proves \eqref{eq:ChebDiffI}. The case of $\xi=\pm1$ can be
proved in a similar manner. This completes the proof.
\end{proof}

\begin{remark}
From Theorem \ref{thm:ChebDiff} we see that the maximum error of
Chebyshev spectral differentiations is actually determined by the
error at one of the endpoints whenever $\xi\in(-1,1)$ and by the
error at the singularity whenever $\xi=\pm1$.
\end{remark}

\begin{remark}
In the case of $\xi\in(-1,1)$, it follows from \eqref{eq:ChebDiffI}
and \eqref{eq:AsyCaseII} that
\[
|f(\xi)-f_n(\xi)|=O(n^{-\alpha}), \quad
|f'(\xi)-f_n{'}(\xi)|=O(n^{-\alpha}).
\]
Interestingly, we see that differentiation does not lead to a
deterioration in the rate of convergence of $f_n$ at the
singularity. Indeed, the estimates of $f-f_n$ and $f{'}-f_n{'}$ as
$n\rightarrow\infty$ are determined by the asymptotics of
$\Psi_{\alpha}^{\mathrm{C}}(2\varphi_{\xi}^{+},n)+\Psi_{\alpha}^{\mathrm{C}}(2\varphi_{\xi}^{-},n)$
and
$\Psi_{\alpha-1}^{\mathrm{S}}(2\varphi_{\xi}^{+},n)+\Psi_{\alpha-1}^{\mathrm{S}}(2\varphi_{\xi}^{-},n)$,
respectively. By Lemma \ref{lem:DecaySeries}, we find that
$\Psi_{\alpha}^{\mathrm{C}}(2\varphi_{\xi}^{+},n)=O(n^{-\alpha-1})$
and
$\Psi_{\alpha-1}^{\mathrm{S}}(2\varphi_{\xi}^{+},n)=O(n^{-\alpha})$
for all $x\in\Omega$, and
$\Psi_{\alpha}^{\mathrm{C}}(2\varphi_{\xi}^{-},n)=O(n^{-\alpha-1})$
and
$\Psi_{\alpha-1}^{\mathrm{S}}(2\varphi_{\xi}^{-},n)=O(n^{-\alpha})$
when $x\in\Omega\setminus\{\xi\}$ and
$\Psi_{\alpha}^{\mathrm{C}}(2\varphi_{\xi}^{-},n)=O(n^{-\alpha})$
and $\Psi_{\alpha-1}^{\mathrm{S}}(2\varphi_{\xi}^{-},n)=0$ when
$x=\xi$. Hence, the asymptotic of
$\Psi_{\alpha}^{\mathrm{C}}(2\varphi_{\xi}^{-},n)$ deteriorates at
$x=\xi$, but this is not the case for
$\Psi_{\alpha-1}^{\mathrm{S}}(2\varphi_{\xi}^{-},n)$.
\end{remark}

\begin{corollary}\label{col:ChebDiffMax}
Let the assumption in Theorem \ref{thm:ChebDiff} be fulfilled. As
$n\rightarrow\infty$, the error estimates of Chebyshev spectral
differentiation in the maximum norm are
\begin{align}\label{eq:ChebDiffMax}
\|f{'}-f_n{'}\|_{L^{\infty}(\Omega)} = \left\{
            \begin{array}{ll}
{\displaystyle O(n^{1-\alpha}) },    & \hbox{$\xi\in(-1,1)$,}   \\[8pt]
{\displaystyle O(n^{2-2\alpha}) },   & \hbox{$\xi=\pm1$.}
            \end{array}
            \right.
\end{align}
Moreover, the maximum error of Chebyshev spectral differentiation is
attained at one of the endpoints whenever $\xi\in(-1,1)$ and at the
singularity $\xi$ whenever $\xi=\pm1$.
\end{corollary}
\begin{proof}
The corollary is a direct consequence of Theorem \ref{thm:ChebDiff}.
\end{proof}

\begin{figure}[ht]
\centering
\includegraphics[width=0.5\textwidth]{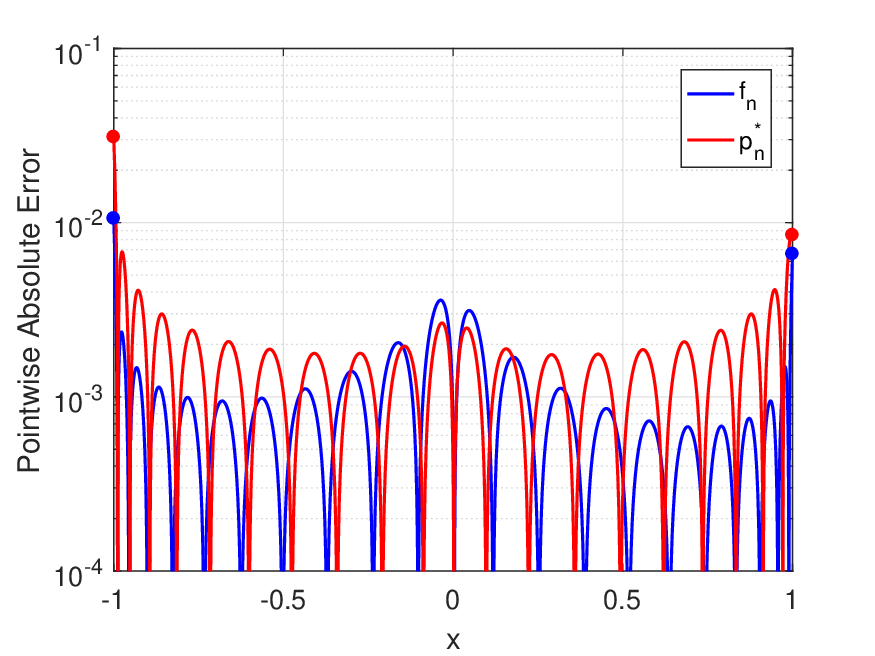}~~
\includegraphics[width=0.5\textwidth]{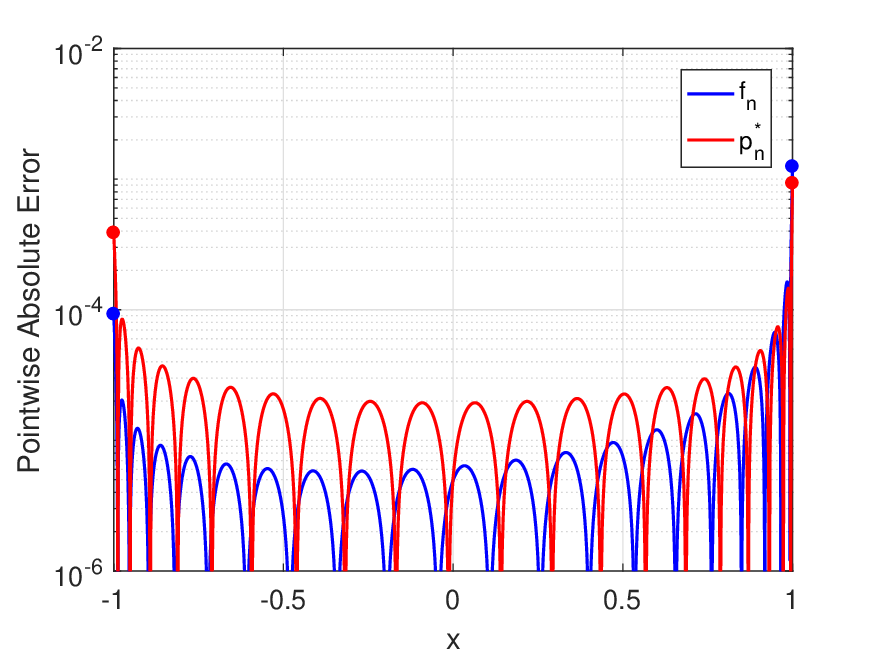}\\
\includegraphics[width=0.5\textwidth]{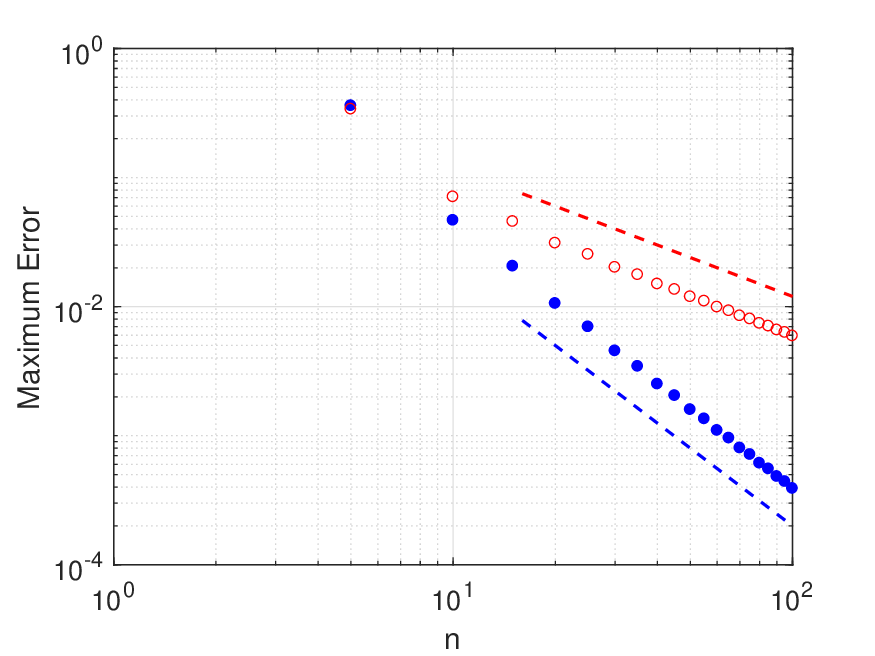}~~
\includegraphics[width=0.5\textwidth]{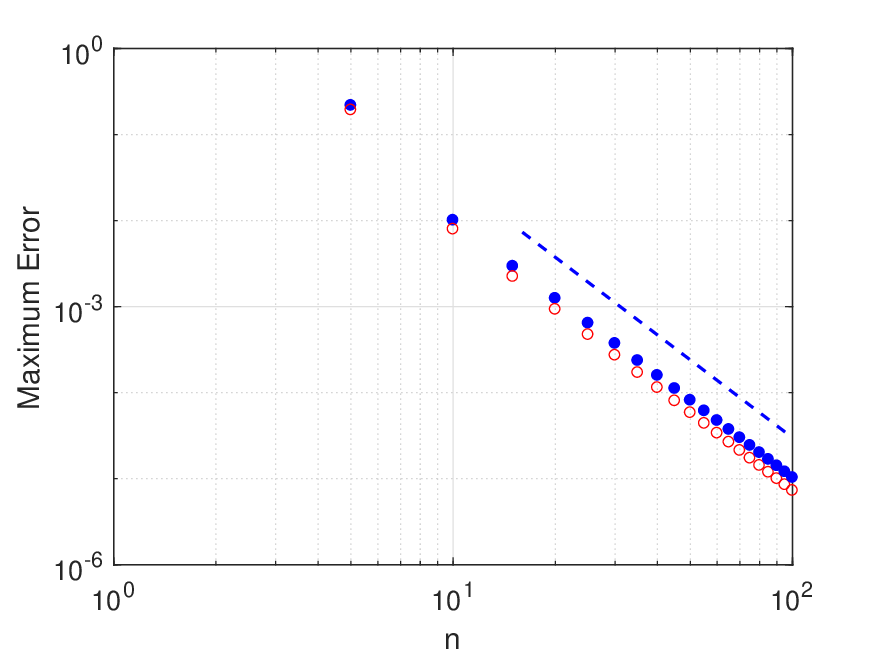}
\caption{Top row shows $|f{'}-f_n{'}|$ and $|f{'}-p_n^{*}{'}|$ for
$f(x)=|x|^{3}e^{\sin(x)}$ (left) and $f(x)=(1-x)^{5/2}e^{\sin(x)}$
(right) with $n=20$, and these points indicate the errors of
$f_n{'}$ and $p_n^{*}{'}$ at both endpoints. Bottom row shows
$\|f{'}-f_n{'}\|_{L^{\infty}(\Omega)}$ (dots) and
$\|f{'}-p_n^{*}{'}\|_{L^{\infty}(\Omega)}$ (circles) for
$f(x)=|x|^{3}e^{\sin(x)}$ (left) and $f(x)=(1-x)^{5/2}e^{\sin(x)}$
(right). The dash lines in the left panel show $O(n^{2-\alpha})$
(upper) and $O(n^{1-\alpha})$ (bottom) and the dash line in the
right panel shows $O(n^{2-2\alpha})$.} \label{fig:ChebDiff1}
\end{figure}

We close this subsection by making a comparison of spectral
differentiations using $f_n$ and $p_n^{*}$. Regarding the maximum
error and pointwise error of spectral differentiations using
$p_n^{*}$, it is still unclear whether they can be found in the
literature. Numerical experiments show that
\begin{align}\label{eq:BestDiffMax}
\|f{'} - p_n^{*}{'} \|_{L^{\infty}(\Omega)} = \left\{
            \begin{array}{ll}
{\displaystyle O(n^{2-\alpha}) },    & \hbox{$\xi\in(-1,1)$,}   \\[8pt]
{\displaystyle O(n^{2-2\alpha}) },   & \hbox{$\xi=\pm1$,}
            \end{array}
            \right.
\end{align}
and the maximum error is always attained at one of the endpoints for
large $n$ (see Figure \ref{fig:ChebDiff1}). With Corollary
\ref{col:ChebDiffMax} and \eqref{eq:BestDiffMax}, we conclude that
the rate of convergence of $f_n$ in the maximum norm is actually one
power of $n$ faster than that of $p_n^{*}$ whenever $\xi\in(-1,1)$
and both are the same whenever $\xi=\pm1$. In Figure
\ref{fig:ChebDiff1} we illustrate the maximum errors and the
pointwise errors of spectral differentiations using $f_n$ and
$p_n^{*}$. As expected, we see that the maximum errors of $f_n$
converges one power of $n$ faster than that of $p_n^{*}$ for
$\xi\in(-1,1)$ and both converge at the same rate whenever
$\xi=\pm1$. Moreover, concerning the pointwise errors of spectral
differentiation using $f_n$ and $p_n^{*}$, we see that Chebyshev
spectral differentiations converge also faster than their best
counterparts except in a small neighborhood of the singularity.

\subsection{Other orthogonal projections}
The error estimate of other orthogonal projections such as Legendre,
Gegenbauer and Jacobi projections has attracted attention in recent
years (see, e.g.,
\cite{babuska2019error,liu2020legendre,wang2020legendre,wang2021gegenbauer,xiang2020jacobi,xiang2021super}).
In \cite[Figure~3]{wang2020legendre} it has been observed that the
error curves of Legendre projections illustrate similar character as
that of Chebyshev projections. In the case of Gegenbauer and Jacobi
projections, however, the situation is slightly more complicated
since the errors at both endpoints may dominate the maximum error
(see \cite{wang2021gegenbauer}). In the following, we restrict our
attention to the case of Legendre projections and show the intrinsic
connection between the error of Legendre projections and the two
functions $\Psi_{\nu}^{\mathrm{C}}(x,n)$ and
$\Psi_{\nu}^{\mathrm{S}}(x,n)$ defined in \eqref{def:PsiFun}.

Let $f$ be the model function defined in \eqref{def:Model} and let
$\mathcal{P}_n(x)$ be its Legendre projection of degree $n$, i.e.,
\begin{align}
\mathcal{P}_n(x) = \sum_{k=0}^{n} a_k^{L} P_k(x), \quad a_k^{L} =
\frac{2k+1}{2} \int_{-1}^{1} f(x) P_k(x) \mathrm{d}x,
\end{align}
where $P_k(x)$ is the Legendre polynomial of degree $k$. By taking
the Taylor series of $g(x)$ at $x=\xi$ and using
\cite[Equation~(7.232.3)]{gradshteyn2007table}, we obtain the
asymptotic behavior of the Legendre coefficients
\begin{align}\label{eq:LegCoeffAsym}
a_k^{L} = \mathcal{E}(\alpha,\xi) \frac{\lambda_1(\xi) T_k(\xi) +
\lambda_2(\xi) (1-\xi^2)^{1/2} U_{k-1}(\xi) }{k^{\alpha+1/2}} +
O(k^{-\alpha-3/2}),
\end{align}
where $\lambda_j(\xi) = (1+\xi)^{1/2} + (-1)^{j+1} (1-\xi)^{1/2}$
and
\begin{align}
\mathcal{E}(\alpha,\xi) = -\sqrt{\frac{2}{\pi}}
(1-\xi^2)^{\alpha/2+1/4} \sin\left(\frac{\alpha\pi}{2}\right)
\Gamma(\alpha+1) g(\xi).  \nonumber
\end{align}
On the other hand, recall from \cite[Theorem~8.21.2]{szego1975orth}
that
\begin{align}\label{eq:LegPolyAsym}
P_k(x) &= \sqrt{\frac{2}{\pi}} \frac{(1-x^2)^{-1/4}}{k^{1/2}}
\cos\left(\frac{2k+1}{2}\arccos(x) - \frac{\pi}{4}
\right) + O(k^{-3/2}) \nonumber \\
&= \frac{(1-x^2)^{-1/4}}{\sqrt{2\pi k}}\left[ \lambda_1(x) T_k(x) +
\lambda_2(x)(1-x^2)^{1/2}U_{k-1}(x) \right] + O(k^{-3/2}),
\end{align}
where the asymptotic formula on the right-hand side of
\eqref{eq:LegPolyAsym} holds uniformly on the interval
$[-1+\varepsilon,1-\varepsilon]$, where $\varepsilon>0$ is small.
Let us now consider the pointwise error estimate of
$\mathcal{P}_n(x)$ and we restrict the analysis to the case where
$\xi\in(-1,1)$ since the case $\xi=\pm1$ can be treated in a similar
way. Combining \eqref{eq:LegCoeffAsym} and \eqref{eq:LegPolyAsym}
and after some elementary calculations, we obtain for $x\in(-1,1)$
that
\begin{align}\label{eq:LegProjError}
f(x) - \mathcal{P}_n(x) % &= \sum_{k=n+1}^{\infty} a_k^L P_k(x) \nonumber \\
= \frac{\mathcal{E}(\alpha,\xi)}{\sqrt{2\pi} (1-x^2)^{1/4}} &\left[
\frac{\lambda_1(\xi) \lambda_1(x)}{2} \sum_{k=n+1}^{\infty} \frac{
\cos(2k\varphi_{\xi}^{+}) + \cos(2k\varphi_{\xi}^{-})}{k^{\alpha+1}}  \right. \nonumber \\
& + \frac{\lambda_1(\xi) \lambda_2(x)}{2} \sum_{k=n+1}^{\infty}
\frac{\sin(2k\varphi_{\xi}^{+}) + \sin(2k\varphi_{\xi}^{-}) }{k^{\alpha+1}} \nonumber \\
& + \frac{\lambda_1(x)\lambda_2(\xi)}{2} \sum_{k=n+1}^{\infty}
\frac{\sin(2k\varphi_{\xi}^{+}) - \sin(2k\varphi_{\xi}^{-}) }{k^{\alpha+1}}  \\
& \left. + \frac{\lambda_2(x)\lambda_2(\xi)}{2}
\sum_{k=n+1}^{\infty} \frac{\cos(2k\varphi_{\xi}^{+}) -
\cos(2k\varphi_{\xi}^{-})}{k^{\alpha+1}} \right] \nonumber \\
& + O(n^{-\alpha-2}). \nonumber
\end{align}
Clearly, we see that the leading term of the error of
$\mathcal{P}_n(x)$ can be represented as a linear combination of
$\Psi_{\nu}^{\mathrm{C}}(x,n)$ and $\Psi_{\nu}^{\mathrm{S}}(x,n)$.
Furthermore, applying Lemma \ref{lem:DecaySeries} to
\eqref{eq:LegProjError}, we can deduce that the term inside the
square bracket on the right-hand side of the last equality behaves
like $O(n^{-\alpha-1})$ whenever $x$ is away from $\xi$ and behaves
like $O(n^{-\alpha})$ whenever $x=\xi$. For $x=\pm1$, combining the
fact that $P_k(\pm1)=(\pm1)^k$ for all $k\in\mathbb{N}_0$ as well as
Lemma \ref{lem:DecaySeries} and \eqref{eq:LegCoeffAsym} yields
\begin{align}
f(\pm1) - \mathcal{P}_n(\pm1) &= \sum_{k=n+1}^{\infty} (\pm1)^k
a_k^L = \mathcal{E}(\alpha,\xi) \left( \lambda_1(\xi)
\sum_{k=n+1}^{\infty} \frac{(\pm1)^k T_k(\xi)}{k^{\alpha+1/2}}
\right. \nonumber \\
&~~~~~~~~~~~ \left.  + \lambda_2(\xi) \sum_{k=n+1}^{\infty}
\frac{(\pm1)^k (1-\xi^2)^{1/2} U_{k-1}(\xi)}{k^{\alpha+1/2}} \right)
+ \cdots.
\end{align}
We see that the leading term on the right-hand side of the last
equality can also be represented by a linear combination of
$\Psi_{\nu}^{\mathrm{C}}(x,n)$ and $\Psi_{\nu}^{\mathrm{S}}(x,n)$.
Combining this with Lemma \ref{lem:DecaySeries} we easily obtain
$f(\pm1) - \mathcal{P}_n(\pm1)=O(n^{-\alpha-1/2})$. Hence, we
conclude that the pointwise error estimate of $\mathcal{P}_n(x)$ is
\begin{align}\label{eq:Leg}
|f(x) - \mathcal{P}_n(x)| = \left\{
            \begin{array}{ll}
{\displaystyle O(n^{-\alpha}) },     & \hbox{if~$x=\xi$,}   \\[8pt]
{\displaystyle O(n^{-\alpha-1/2}) }, & \hbox{if~$x=\pm1$,}   \\[8pt]
{\displaystyle O(n^{-\alpha-1}) },   & \hbox{otherwise,}
            \end{array}
            \right.
\end{align}
Comparing the rate of pointwise convergence of $\mathcal{P}_n(x)$
with $p_n^{*}(x)$, we see that $\mathcal{P}_n(x)$ converges faster
than $p_n^{*}(x)$ by a factor of $n$ whenever
$x\in(-1,\xi)\cup(\xi,1)$ and by a factor of $n^{1/2}$ whenever
$x=\pm1$ and both $\mathcal{P}_n(x)$ and $p_n^{*}(x)$ converge at
the same rate whenever $x=\xi$, and these findings justify the error
localization of Legendre projections.

\begin{remark}
In \cite{wang2020legendre} the present author has studied the
optimal rates of convergence of Legendre projections in the maximum
norm for analytic and piecewise analytic functions and function of
the form \eqref{def:Model}. In the case of functions with an
interior singularity, i.e., $\xi\in(-1,1)$, it has been shown that
$\mathcal{P}_n(x)$ and $p_n^{*}(x)$ converge at the same rate based
on the hypothesis that the maximum error of $\mathcal{P}_n(x)$ is
attained at $x=\xi$. Here our pointwise error estimates in
\eqref{eq:Leg} provide a theoretical justification for this
hypothesis.
\end{remark}

\begin{remark}
Note that we have used the uniform asymptotic expansion of Legendre
polynomials in compact subsets of $(-1,1)$ in the derivation of
\eqref{eq:Leg}. A similar approach for the analysis of
polyharmonic-Neumann expansions was used in
\cite[Chapter~3]{adcock2010b}.
\end{remark}

\section{Concluding remarks}\label{sec:conclusion}
We have presented a thorough analysis of the pointwise error
estimate of Chebyshev spectral approximations for functions with a
singularity. Our key finding is that there exists an intrinsic
connection between pointwise error estimates of Chebyshev spectral
approximations and the two functions $\Psi_{\nu}^{\mathrm{C}}(x,n)$
and $\Psi_{\nu}^{\mathrm{S}}(x,n)$ defined in \eqref{def:PsiFun}.
This connection allows us to justify rigorously the error
localization of Chebyshev spectral approximations, such as why their
maximum error is always attained in a small neighborhood of the
singularity and why they converge faster than their best
counterparts except for a small neighborhood of the singularity. We
further extended the framework to Chebyshev spectral
differentiations and Legendre projections and justified their error
localization property using similar arguments.

Finally, we highlight several directions for future research. First,
it is possible to extend the current study to more general forms of
singular functions, such as
\[
f(x)=\sum_{k=1}^{m} |x-\xi_k|^{\alpha_k}g_k(x), \quad  f(x) =
g_0(x)\prod_{k=1}^{m} |x-\xi_k|^{\alpha_k},
\]
where $\{\xi_k\}_{k=1}^{m}\subseteq\Omega$ and
$\{\alpha_k\}_{k=1}^{m}$ are all positive and not even integers.
More precisely, the extension to the former is direct since it is a
linear combination of the model function \eqref{def:Model}. The
extension to the latter is also applicable since it can also be
decomposed into a linear combination of the model function
\eqref{def:Model} by introducing neutralizers (see
\cite[Equation~(4.17)]{tuan1972}). Second, we have shown that
Chebyshev spectral differentiations converge faster than their best
counterparts except in a neighborhood of the singularity and, in the
particular case of functions with an interior singularity, they
converge even faster than their best counterparts in the maximum
norm. Note that spectral differentiations are widely used in the
context of spectral methods, it is of great interest to extend the
study to more general settings, such as Gegenbauer and Jacobi
spectral differentiations.

\section*{Acknowledgements}
This work was supported by National Natural Science Foundation of
China under grant number 11671160. The author would like to thank
Prof. Nick Trefethen for his suggestions that improved the
presentation of the paper. He would also like to thank Prof.
Chengming Huang for several stimulating discussions and the
anonymous referee for helpful comments on this work.

%\appendix

\end{document}